\documentclass[review]{elsarticle}

\usepackage[T1]{fontenc}
\usepackage{amsmath,amssymb,amsfonts}
\usepackage{graphicx}
\usepackage{booktabs}
\usepackage{url}
\usepackage{hyperref}

\newcommand{\rangOpere}{6}          
\newcommand{\nSim}{576}             
\newcommand{\nTrain}{403}           
\newcommand{\nVal}{58}              
\newcommand{\nTest}{115}            
\newcommand{\nNodes}{60914}         
\newcommand{\nParam}{5}             
\newcommand{\PloadMin}{100}   
\newcommand{\PloadMax}{600}   
\newcommand{\MassMin}{50}     
\newcommand{\MassMax}{300}    
\newcommand{\VinitMin}{40}    
\newcommand{\VinitMax}{90}    
\newcommand{\TdiscMin}{20}    
\newcommand{\TdiscMax}{200}   
\newcommand{\TambMin}{-40}    
\newcommand{\TambMax}{50}     
\newcommand{\PloadLev}{4}
\newcommand{\MassLev}{4}
\newcommand{\VinitLev}{4}
\newcommand{\TdiscLev}{3}
\newcommand{\TambLev}{3}
\newcommand{\budTruncMedT}{0.0765}  
\newcommand{\budTruncPToT}{0.2726}  
\newcommand{\budRegMedT}{0.0098}    
\newcommand{\budRegPToT}{0.033}     
\newcommand{\budRealMedT}{0.0752}   
\newcommand{\budRealPToT}{0.2629}   
\newcommand{\budTruncMedS}{0.0197}  
\newcommand{\budTruncPToS}{0.0864}  
\newcommand{\budRegMedS}{0.0005}    
\newcommand{\budRegPToS}{0.0038}    
\newcommand{\budRealMedS}{0.0196}   
\newcommand{\budRealPToS}{0.0858}   
\newcommand{\budDomFracT}{0.86}     
\newcommand{\budDomFracS}{1.0}      
\newcommand{\budSousPlancherFrac}{1.0} 
\newcommand{\tauRel}{0.01}          
\newcommand{\medRelExtTestT}{0.000055} 
\newcommand{\maxRelExtTestT}{0.000438} 
\newcommand{\absExtTestT}{0.3844}   
\newcommand{\conformFracTestT}{1.0} 
\newcommand{\medRelExtTestS}{0.00025}  
\newcommand{\maxRelExtTestS}{0.017417} 
\newcommand{\absExtTestS}{0.5172}   
\newcommand{\conformFracTestS}{0.965}  
\newcommand{\nHorsCritTestS}{4}     
\newcommand{\zoneAireMoyS}{0.3121}  
\newcommand{\zoneAireCVS}{0.4176}   
\newcommand{\zoneRhoAireS}{0.9995}  
\newcommand{\corrHorsDiagS}{0.786}    
\newcommand{\cpLoExtValS}{0.788}      
\newcommand{\cpHiExtValS}{0.961}      
\newcommand{\genNHgPoints}{15}       
\newcommand{\genMedRelHgT}{0.0210}   
\newcommand{\genMaxRelHgT}{0.0744}   
\newcommand{\genConfHgT}{0.2667}     
\newcommand{\genMedRelHgS}{0.0057}   
\newcommand{\genMaxRelHgS}{0.0362}   
\newcommand{\genConfHgS}{0.8000}     
\newcommand{\covExtValT}{1.0}       
\newcommand{\covExtValS}{0.897}     
\newcommand{\femMoyS}{332.5}        
\newcommand{\femCumulS}{160271}     
\newcommand{\tInfMedMs}{2.238}      
\newcommand{\accelMed}{88921}       
\newcommand{\accelCumul}{148579}    
\newcommand{\accelMin}{16980}       
\newcommand{\accelMax}{904399}      
\newcommand{\modeUnPicT}{0.9004}    
\newcommand{\modeUnPicS}{0.9589}    
\newcommand{\rankPicNNNS}{3}        
\newcommand{\epshS}{11.85}          
\newcommand{\epshRelT}{0.002}       
\newcommand{\meshCoarseNodes}{10645}   
\newcommand{\meshProdNodes}{60914}     
\newcommand{\meshFineNodes}{391480}    
\newcommand{\meshCoarseSize}{15}        
\newcommand{\meshProdSize}{8}           
\newcommand{\meshFineSize}{4}           
\newcommand{\svmMaxCoarse}{107.22}  
\newcommand{\svmMaxProd}{213.84}    
\newcommand{\svmMaxFine}{225.70}    
\newcommand{\epshRelTMedian}{0.0002}  
\newcommand{\epshRelTCoin}{0.0023}    
\newcommand{\epshRelSMin}{0.0525}     
\newcommand{\epshRelSMax}{0.0640}     
\newcommand{\ksPvalTZero}{0.0273}   
\newcommand{\ksEcartMoyTZero}{21}   
\newcommand{\ksPlageTZero}{180}     
\newcommand{\ksMinTZero}{20}        
\newcommand{\ksMaxTZero}{200}       
\newcommand{\nomLevel}{0.95}   

\makeatletter
\def\ps@pprintTitle{%
  \let\@oddhead\@empty
  \let\@evenhead\@empty
  \def\@oddfoot{\footnotesize\itshape
    Preprint submitted to Elsevier \hfill July 2026}%
  \let\@evenfoot\@oddfoot}
\makeatother

\usepackage{float}              
\usepackage[section]{placeins}  

\begin{document}

\begin{frontmatter}

\title{Surrogate-to-code verification of a non-intrusive POD-GPR machine-learning emulator of peak thermomechanical fields, with application to a carbon-carbon aircraft brake disc}

\author{Franklin Kamche}
\ead{franklin.kamche@gmail.com}
\affiliation{organization={Independent Researcher},
             city={Paris}, country={France}}

\begin{abstract}
While non-intrusive reduced-order emulators successfully substitute for costly finite-element simulations, certification-grade applications demand a rigorous error assessment relative to the solver's own discretization uncertainty. This work introduces a verification framework designed to evaluate prediction accuracy directly against this numerical uncertainty floor, rather than physical observation. We demonstrate this approach on a carbon-carbon aircraft brake disc during a rejected take-off scenario, predicting peak temperature and peak von Mises stress fields. Operating over a fixed mesh of \nNodes{} nodes, the emulator relies on a pipeline coupling proper orthogonal decomposition (POD) and Gaussian process regression (GPR) parameterized by \nParam{} inputs. We formulate an algebraically exact, a posteriori error budget that isolates the emulator's signed error into distinct truncation and regression components at the extremum. This error is then evaluated against a discretization floor established via mesh-convergence studies on the peak quantities of interest. On a completely independent test set, our emulator matches the high-fidelity solver's peak values within its own discretization uncertainty for both fields, a benchmark we define as 'surrogate-to-code' fidelity. Notably, the residual error stems primarily from linear reduction limits rather than regression inaccuracies. Beyond the main framework, we present two secondary contributions: the concurrent, non-intrusive prediction of the thermal and mechanical peaks, and a probabilistic approach to spatial zone localization for these extrema. This methodology applies to non-intrusive reduced-order models of parameterized finite-element simulations wherever a discretization-error estimate for the quantity of interest can be obtained.
\end{abstract}

\begin{keyword}
surrogate-to-code verification \sep non-intrusive reduced-order model \sep proper orthogonal decomposition \sep Gaussian process regression \sep discretization error \sep carbon-carbon aircraft brake
\end{keyword}

\end{frontmatter}

\section{Introduction}
\label{sec:intro}

Non-intrusive reduced-order emulators substitute for repeated executions of expensive
finite-element solvers by training cheap approximations on solver snapshots. This
approach enables design studies and many-query applications that remain computationally
intractable with the solver alone. When targeted at certification, the utility of such
emulators depends on establishing rigorous error bounds. These bounds require evaluation
against the correct reference. Rather than physical truth, the reference is the
underlying numerical solver. This solver is itself an approximation of physical reality
and possesses its own discretization error. Reducing the emulator error below this
discretization threshold merely refines agreement with the numerical model beyond its own
resolution, which lacks physical meaning.

Intrusive reduced-basis methods benefit from established a posteriori error estimators
that access the governing operators to bound the reduced solution
rigorously~\cite{rozza2008rbaposteriori,grepl2005aposteriori}. In contrast, the
non-intrusive setting operates strictly on input-output data. Because residual-based
estimators are unavailable in this context, establishing a controlled,
certification-grade error assessment against the solver's numerical uncertainty remains
an open challenge. Standard non-intrusive reduced-order models construct maps via
regression of reduced coordinates or via operator
inference~\cite{peherstorfer2016operatorinference,guo2019timedependent,benner2015survey}.
Separate uncertainty and calibration frameworks map the discrepancy between models and
physical systems~\cite{kennedy2001bayesiancalib} within the standard lexicon of
verification and validation~\cite{roy2011vvuq}.

This study evaluates a high-consequence benchmark from aircraft braking. A rejected
take-off converts the entire kinetic energy of the aircraft into thermal energy, causing
the carbon-carbon brake disc to experience a peak temperature followed closely by peak
thermoelastic stress~\cite{ramachandra2026wbreview}. These peak values dictate landing
gear safety and brake reuse. Predicting them across the full operating envelope represents
a primary design and certification requirement. High-fidelity finite-element models of
the disc provide accurate simulations. Their computational cost makes dense parametric
exploration impractical for multi-query studies~\cite{ramachandra2026wbreview}. This
limitation motivates data-driven surrogates for virtual sensing of quantities that are hard to measure directly~\cite{ghienne2023virtualsensors}, prediction of the brake stress state from simulation data~\cite{melgarejo2024virtualstress}, and real-time structural assessment~\cite{mainini2015surrogate}.

This paper presents a verification methodology for non-intrusive
reduced-order emulators derived from high-fidelity solvers, demonstrated on the aircraft
brake disc. The underlying emulator uses a standard pipeline coupling proper orthogonal
decomposition (POD) and Gaussian process regression (GPR). The novelty lies in the
verification framework. We introduce an a posteriori peak error budget that decomposes
the signed emulator error, via an algebraically exact identity, into a reduced-basis
truncation term and a regression term. The observed extremum error is then evaluated
against a discretization floor derived from a mesh-convergence study of the reference
model. Predicting peak values is a difficult task. The extremum is a non-smooth
functional that amplifies local spatial errors, making peak verification more stringent
than global field-wise norms. Our main finding shows that, on a held-out test set, the
emulator reproduces the high-fidelity peak predictions within the solver's own
discretization uncertainty. We define this level of agreement as surrogate-to-code
fidelity, which extends code verification concepts to surrogate models rather than
asserting physical validation. Two secondary contributions support this framework. We
perform concurrent non-intrusive predictions of thermal and mechanical peaks within a
single workflow. We also introduce a probabilistic localization of spatial extrema zones,
propagating Gaussian process posterior uncertainties directly onto the computational mesh, reported under an explicit appraisal of its added value.

The manuscript is organized as follows. Section~\ref{sec:problem} outlines the
mathematical formulation and the peak quantities of interest. Section~\ref{sec:hifi}
defines the high-fidelity model, the dataset, and the discretization floor.
Section~\ref{sec:emulator} details the non-intrusive emulator pipeline.
Section~\ref{sec:budget} derives the peak error budget and the level-one uncertainty
metrics. Section~\ref{sec:results} presents the results obtained on the held-out test
set. Section~\ref{sec:ablation} provides an ablation study, evaluates the scope of the linear reduction, and discusses model interpretability. Section~\ref{sec:discussion}
addresses limitations and positioning, and Section~\ref{sec:conclusion} concludes the
paper.
\section{Problem statement and quantities of interest at peak}
\label{sec:problem}

Brake certification protocols rely on the peak thermal state of a brake disc to evaluate structural integrity. This maximum thermal state also guides operational decisions regarding brake reuse and aircraft turnaround times after an aborted take-off. A rejected take-off represents the most severe thermal loading case for sizing braking systems. Modern high-energy aircraft designs utilize carbon-carbon composite discs to exploit their high thermal capacity and minimal mass~\cite{ramachandra2026wbreview}. In these scenarios, the brakes convert the kinetic energy of the aircraft entirely into thermal energy, and this heat generation drives the subsequent mechanical response of the disc. This study focuses on this peak operating state. We develop a non-intrusive emulator to predict both the temperature field and the von Mises stress field of a carbon-carbon brake disc at their respective peaks using only the operating parameters as inputs. The model also provides predictive uncertainty estimates for these peak quantities relative to the high-fidelity solver predictions.

\begin{figure}[htbp]
  \centering
  \includegraphics[width=\linewidth]{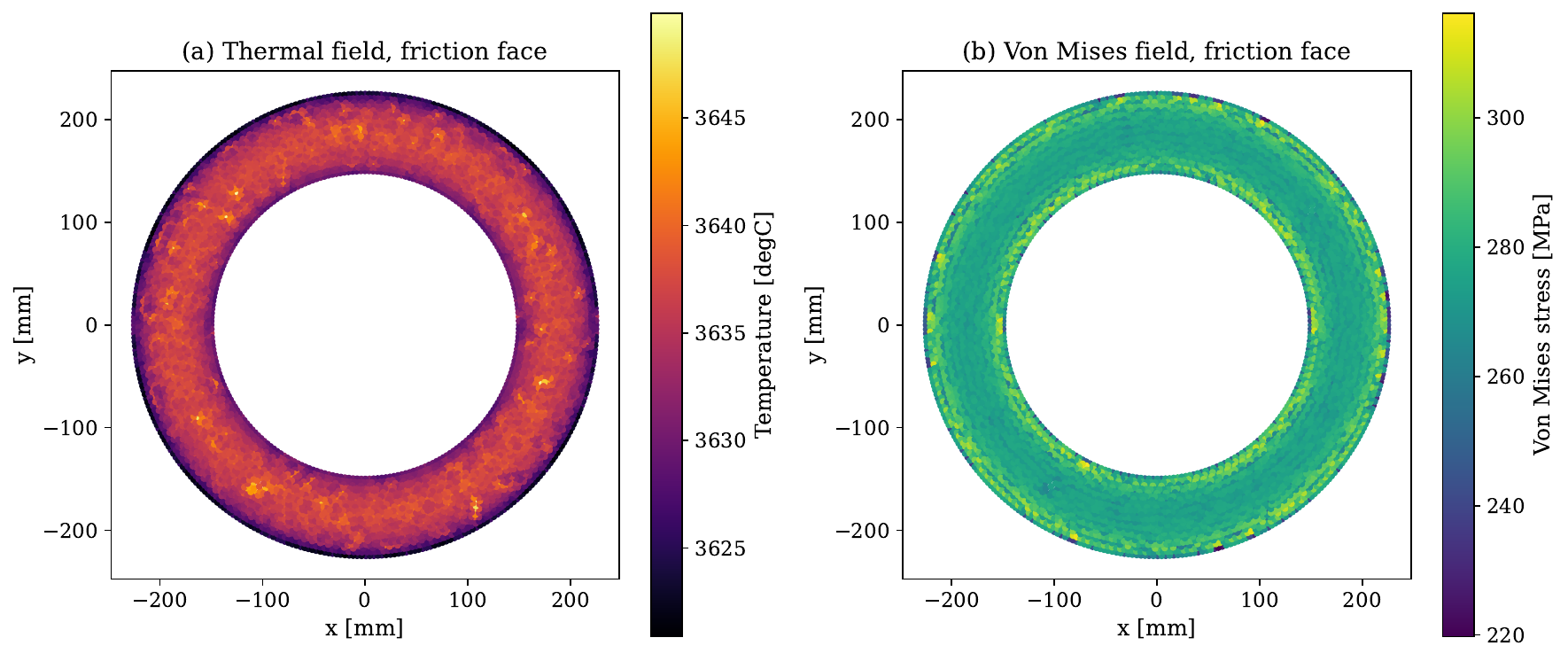}
  \caption{Spatial distribution of field extrema across the disc friction surface. The maximum values for the peak temperature and peak von Mises stress fields concentrate within a nearly flat annular band. Although the peak value remains a stable quantity of interest, the specific node where this maximum occurs shifts within the band, making single-node identification unreliable. Extrema localization is therefore defined as a spatial zone rather than a single node.}
  \label{fig:extremum-topology}
\end{figure}

\subsection{Parametric input space}
\label{subsec:inputs}

We define each operating condition by a vector $\mathbf{x} \in \mathbb{R}^{\nParam}$ containing \nParam{} scalar parameters:
\begin{equation}
  \mathbf{x} = \bigl(P,\; M,\; V_0,\; T_0,\; T_{\mathrm{amb}}\bigr),
  \label{eq:input-vector}
\end{equation}
where $P$ represents the braking load in kilonewtons, $M$ is the aircraft mass in tonnes, $V_0$ is the initial speed during the rejected take-off in meters per second, $T_0$ denotes the initial disc temperature in degrees Celsius, and $T_{\mathrm{amb}}$ is the ambient air temperature in degrees Celsius. Section~\ref{sec:hifi} details the admissible parameter ranges and the sampling strategy used to build the dataset.

\subsection{Prediction targets}
\label{subsec:targets}

For each operating point, the high-fidelity model described in Section~\ref{sec:hifi} computes the transient thermal and mechanical response over a fixed mesh containing $N_{\mathrm{n}} = \nNodes{}$ nodes. This shared discretization guarantees direct nodal correspondence across all dataset samples. We denote the transient temperature and von Mises stress fields as $\mathcal{T}(\mathbf{x},t) \in \mathbb{R}^{N_{\mathrm{n}}}$ and $\mathcal{S}(\mathbf{x},t) \in \mathbb{R}^{N_{\mathrm{n}}}$. Rather than emulating the complete transient history, we extract each field at the specific time step where its global spatial maximum is reached:
\begin{equation}
  t_T(\mathbf{x}) = \arg\max_{t}\ \max_{1 \le i \le N_{\mathrm{n}}}
    \mathcal{T}_i(\mathbf{x},t),
  \qquad
  t_\sigma(\mathbf{x}) = \arg\max_{t}\ \max_{1 \le i \le N_{\mathrm{n}}}
    \mathcal{S}_i(\mathbf{x},t).
  \label{eq:peak-instants}
\end{equation}
Focusing on these discrete instances aligns with certification requirements, which prioritize peak thermal and mechanical loads over full transient histories. The resulting prediction targets are the spatial fields:
\begin{equation}
  T(\mathbf{x}) = \mathcal{T}\bigl(\mathbf{x}, t_T(\mathbf{x})\bigr)
    \in \mathbb{R}^{N_{\mathrm{n}}},
  \qquad
  \sigma(\mathbf{x}) = \mathcal{S}\bigl(\mathbf{x}, t_\sigma(\mathbf{x})\bigr)
    \in \mathbb{R}^{N_{\mathrm{n}}},
  \label{eq:target-fields}
\end{equation}
representing the nodal temperature (in degrees Celsius) and the von Mises stress (in megapascals) at their respective peak times. In our dataset, the mechanical peak consistently occurs at or after the thermal peak, yielding $t_\sigma(\mathbf{x}) \ge t_T(\mathbf{x})$. This is an empirical property of the present dataset and is not imposed by the formulation. The emulator learns these two static states directly, bypassing trajectory modeling.

\subsection{Certified quantities of interest}
\label{subsec:qoi}

The target certified values are scalar output functionals derived from the predicted spatial fields, following established principles for parametrized problems~\cite{prudhomme2002outputbound}. For any nodal field $u \in \mathbb{R}^{N_{\mathrm{n}}}$, we define the peak operator as:
\begin{equation}
  Q(u) = \max_{1 \le i \le N_{\mathrm{n}}} u_i ,
  \label{eq:peak-functional}
\end{equation}
yielding the two primary quantities of interest, $Q\bigl(T(\mathbf{x})\bigr)$ and $Q\bigl(\sigma(\mathbf{x})\bigr)$. This functional $Q$ represents a non-differentiable, nonlinear $L^{\infty}$ extremum evaluated over the mesh nodes. Operating on discrete vectors makes $Q$ a mesh-dependent approximation of the continuous spatial supremum, supporting our focus on a surrogate-to-code verification strategy. Standard reduced-basis output bounding focuses on linear functionals within coercive systems~\cite{prudhomme2002outputbound}. Since $Q$ is a nonlinear maximum operator, we certify the peak values indirectly using a posteriori bounds on the spatial field error instead of direct analytical output bounds. The high-fidelity finite-element simulations in Section~\ref{sec:hifi} serve as the reference. Verifying the emulator directly against solver outputs establishes a surrogate-to-code rather than a surrogate-to-physics certification. Section~\ref{sec:budget} develops the exact a posteriori peak error budget supporting this verification.

\subsection{Spatial structure of the extremum}
\label{subsec:topology}

While the peak scalar $Q(u)$ remains mathematically well-posed, the specific node where the maximum occurs lacks robustness as a physical descriptor. The high-fidelity fields in Figure~\ref{fig:extremum-topology} show that the extrema occupy a broad, nearly flat annular band on the disc friction surface. The maximum value varies smoothly across the operating space, whereas the argmax coordinates shift abruptly within this flat crest, preventing reliable node-level identification. This flat profile near the maximum justifies identifying the extremum location as a spatial region rather than a single nodal point. For any field $u$ and tolerance parameter $\alpha \in (0,1)$, we define this region as the super-level set:
\begin{equation}
  Z_\alpha(u) = \bigl\{\, i \in \{1,\dots,N_{\mathrm{n}}\}
    \;\mid\; u_i \ge (1-\alpha)\,Q(u) \,\bigr\}.
  \label{eq:zone}
\end{equation}
We develop the probabilistic prediction of this spatial zone as the secondary contribution of this study, with an explicit caveat on its added value in Section~\ref{sec:ablation}. The formal a posteriori peak error budget derived in Section~\ref{sec:budget} constitutes the primary methodology of our verification framework.
\section{High-fidelity model, dataset, and the discretization floor}
\label{sec:hifi}

This section details the reference finite-element model used to train and test the emulator. We outline the dataset generation, the specific modeling assumptions bounding the surrogate-to-code claim, and the discretization floor that defines the limit of meaningful emulator accuracy.

\subsection{Governing thermomechanical model}
\label{subsec:governing}

A rejected take-off converts the kinetic energy of the aircraft into heat at the disc friction interface, generating a transient temperature field that subsequently drives thermoelastic stress. This coupling is unidirectional, meaning temperature changes alter the stress state while mechanical work does not affect the thermal field. The disc occupies a three-dimensional domain $\Omega \subset \mathbb{R}^{3}$.

The thermal field satisfies the transient heat equation. With density $\rho$, specific heat $c_p$, and orthotropic conductivity tensor $\mathbf{k}$, the energy balance and Fourier's law yield
\begin{equation}
  \rho\, c_p\, \frac{\partial T}{\partial t}
    = \nabla \cdot \bigl( \mathbf{k}\, \nabla T \bigr).
  \label{eq:heat}
\end{equation}
The carbon-carbon composite is orthotropic, with a diagonal conductivity tensor in the material axes, $\mathbf{k} = \operatorname{diag}(k_1, k_2, k_3)$, reflecting differing in-plane and through-thickness transport. Frictional heating acts strictly as a boundary flux, meaning no volumetric source term enters \eqref{eq:heat} and no mechanical dissipation is added. We assume temperature-independent material properties over the analyzed operating range. The initial condition imposes a uniform temperature, $T = T_0$ at $t = 0$. Boundary conditions combine the frictional heat flux on contact surfaces with convective cooling on free surfaces,
\begin{equation}
  \begin{aligned}
    -\mathbf{k}\, \nabla T \cdot \mathbf{n} &= \dot{q}_0(t)
      && \text{(Neumann, contact faces)}, \\
    -\mathbf{k}\, \nabla T \cdot \mathbf{n} &= h\,(T - T_{\mathrm{amb}})
      && \text{(Robin, free faces)}.
  \end{aligned}
  \label{eq:thermal-bc}
\end{equation}
The transient flux magnitude $\dot{q}_0(t)$ scales with the dissipated kinetic energy, decreasing linearly to zero over the braking duration, while the convective term references the ambient temperature $T_{\mathrm{amb}}$.

Quasi-static equilibrium dictates the stress field, $\nabla \cdot \boldsymbol{\sigma} = \mathbf{0}$, neglecting inertial terms. We close the system using an orthotropic thermoelastic constitutive law. The thermal expansion tensor $\boldsymbol{\alpha} = \operatorname{diag}(\alpha_1, \alpha_2, \alpha_3)$ is also aligned with the material axes, yielding
\begin{equation}
  \boldsymbol{\sigma}
    = \mathbf{C} : \bigl( \boldsymbol{\varepsilon}
      - \boldsymbol{\alpha}\,(T - T_{\mathrm{ref}}) \bigr),
  \label{eq:thermoelastic}
\end{equation}
where $\mathbf{C}$ is the orthotropic stiffness tensor and $T_{\mathrm{ref}}$ is the stress-free reference temperature. This reference value differs from the initial temperature if the disc contains residual thermal stress at start-up. We focus on the scalar von Mises equivalent stress,
\begin{equation}
  \sigma_{\mathrm{vm}} = \sqrt{\tfrac{3}{2}\, \mathbf{s} : \mathbf{s}},
  \qquad
  \mathbf{s} = \boldsymbol{\sigma}
    - \tfrac{1}{3}\operatorname{tr}(\boldsymbol{\sigma})\, \mathbf{I}.
  \label{eq:vonmises}
\end{equation}
Because the stress response lags the thermal load, the two fields reach their maxima at distinct times, as defined by the peak operators in \eqref{eq:peak-instants}. The repository accompanying this study contains the complete set of orthotropic material properties.

\subsection{Finite-element reference and solver}
\label{subsec:fem}

We solve the weak Galerkin form of equations \eqref{eq:heat} to \eqref{eq:vonmises} on a fixed three-dimensional mesh of $N_{\mathrm{n}} = \nNodes{}$ nodes. The simulation uses CalculiX~\cite{dhondt2004calculix} to perform an implicit, fully coupled temperature-displacement analysis, resolving quasi-static mechanical equilibrium via the PaStiX sparse direct solver~\cite{henon2002pastix}. Element definitions, time-stepping parameters, and convergence tolerances match the deposited input decks. For any input vector $\mathbf{x}$, the solver outputs transient history fields $\mathcal{T}(\mathbf{x},t)$ and $\mathcal{S}(\mathbf{x},t)$, which yield the peak fields $T(\mathbf{x})$ and $\sigma(\mathbf{x})$ of Section~\ref{sec:problem}. We compare all emulator outputs directly against this reference solution using the identical mesh layout.

The high computational cost of individual finite-element runs prevents dense, direct parametric sweeps of the input space, justifying a data-driven surrogate approach. We quantify this computational acceleration in Section~\ref{sec:results}. The simulation cost scales with the number of transient increments, which scales with the total braking duration. Higher braking loads shorten this duration, meaning simulations at high braking load require less computational time.

\subsection{Design of experiments and dataset}
\label{subsec:doe}

\begin{table}[tbp]
  \centering
  \caption{Operational parameters and full factorial design. The ranges and the level counts define the sampled input space.}
  \label{tab:doe}
  \begin{tabular}{llccl}
    \toprule
    Parameter & Symbol & Unit & Range & Levels \\
    \midrule
    Braking load        & $P$                & kN          & \PloadMin{} to \PloadMax{} & \PloadLev{} \\
    Aircraft mass       & $M$                & t           & \MassMin{} to \MassMax{}   & \MassLev{}  \\
    Initial speed       & $V_0$              & m/s         & \VinitMin{} to \VinitMax{} & \VinitLev{} \\
    Initial temperature & $T_0$              & $^{\circ}$C & \TdiscMin{} to \TdiscMax{} & \TdiscLev{} \\
    Ambient temperature & $T_{\mathrm{amb}}$ & $^{\circ}$C & \TambMin{} to \TambMax{}   & \TambLev{}  \\
    \bottomrule
  \end{tabular}
\end{table}

We sample the parametric input space using a full factorial design~\cite{santner2003doe}. Table~\ref{tab:doe} outlines the operational ranges and level allocations for each input parameter. This factorial grid defines the boundaries for model training and verification. While computationally demanding to construct initially, the structured grid provides uniform coverage of the parameter space and captures multi-parameter coupling. Crossing \PloadLev{} loads, \MassLev{} masses, \VinitLev{} initial velocities, \TdiscLev{} initial temperatures, and \TambLev{} ambient temperatures yields a dataset of $\nSim{}$ high-fidelity simulations. Generating all cases on the identical production mesh guarantees exact node-to-node correspondence.

The transient history records starting from the first saved increment, excluding the initial state at $t = 0$. To account for variable braking durations across cases, we map the fields to a uniform time grid and mask inactive increments past the stopping time. The peak operators in \eqref{eq:peak-instants} evaluate only this active, unmasked window.

A fixed random partition divides the $\nSim{}$ cases into $\nTrain{}$ training, $\nVal{}$ validation, and $\nTest{}$ testing samples. We establish this split prior to model development and preserve it throughout the study. The test set is opened a single time, for the final metrics of Section~\ref{sec:results}, and is never used for tuning or model selection.

The dataset carries a documented representativity limit. Because the initial temperature $T_0$ is sampled at discrete factorial levels, we use a two-sample Kolmogorov-Smirnov comparison of its empirical distribution between the test and training splits as an approximate distributional diagnostic rather than an exact hypothesis test. The p-value of $\ksPvalTZero{}$ signals a mild distributional shift, with a mean discrepancy of $\ksEcartMoyTZero{}$ degrees over the $\ksPlageTZero{}$ degree range (from $\ksMinTZero{}$ to $\ksMaxTZero{}$ degrees). Both splits share identical boundary limits. We treat this statistical shift as a known validation limit and analyze its impact in Section~\ref{sec:results}.

\subsection{Limitations of the high-fidelity model and surrogate-to-code scope}
\label{subsec:limits}

The reference model approximates the qualified engineering setup rather than replicating the physical brake disc. Frictional heating enters the system as a prescribed boundary flux with a predetermined temporal profile, bypassing a fully coupled thermomechanical contact formulation that would dynamically link local pressure to heat generation~\cite{wriggers2006contact}. Additionally, we assume temperature-independent material properties and neglect radiative heat transfer relative to convection. Consequently, the emulator cannot predict behavior under altered cooling environments or ventilation designs where hot-spot locations might shift~\cite{ramachandra2026wbreview}.

These modeling assumptions restrict the physical scope of the reference. The emulator is verified directly against the finite-element solver outputs, establishing a surrogate-to-code agreement. The error budget in Section~\ref{sec:budget} measures only the distance between the surrogate and the numerical reference, leaving the physical gap to the real brake disc outside the scope of this study.

\subsection{Discretization floor}
\label{subsec:floor}

The finite-element reference contains inherent numerical discretization error, representing a baseline uncertainty in the sense of verification and validation~\cite{roy2011vvuq}. This discretization error defines a threshold below which further reduction in emulator error lacks physical or numerical meaning. We quantify this floor using three nested meshes with characteristic element sizes of $\meshCoarseSize{}$, $\meshProdSize{}$, and $\meshFineSize{}$ millimeters, corresponding to $\meshCoarseNodes{}$, $\meshProdNodes{}$, and $\meshFineNodes{}$ nodes. The peak von Mises stress increases from $\svmMaxCoarse{}$ to $\svmMaxProd{}$ and then to $\svmMaxFine{}$ megapascals, demonstrating monotonic convergence where the coarsest mesh under-resolves the stress concentrations. Because the coarse mesh lies outside the asymptotic convergence regime, we avoid Richardson extrapolation and estimate the discretization floor using only the two finer meshes. For the peak temperature, the relative difference between the production and refined meshes is $\epshRelT{}$, which is negligible compared to the emulator errors in Section~\ref{sec:results}. For the peak von Mises stress, the relative difference is more sensitive and case-dependent, ranging from $\epshRelSMin{}$ to $\epshRelSMax{}$ (equal to $\epshS{}$ megapascals at the reference case) due to localized spatial gradients. We incorporate this stress-specific discretization floor as an explicit threshold in the error budget derived in Section~\ref{sec:budget}.
\section{Non-intrusive POD-GPR emulator}
\label{sec:emulator}

\begin{figure}[htbp]
  \centering
  \includegraphics[width=\linewidth]{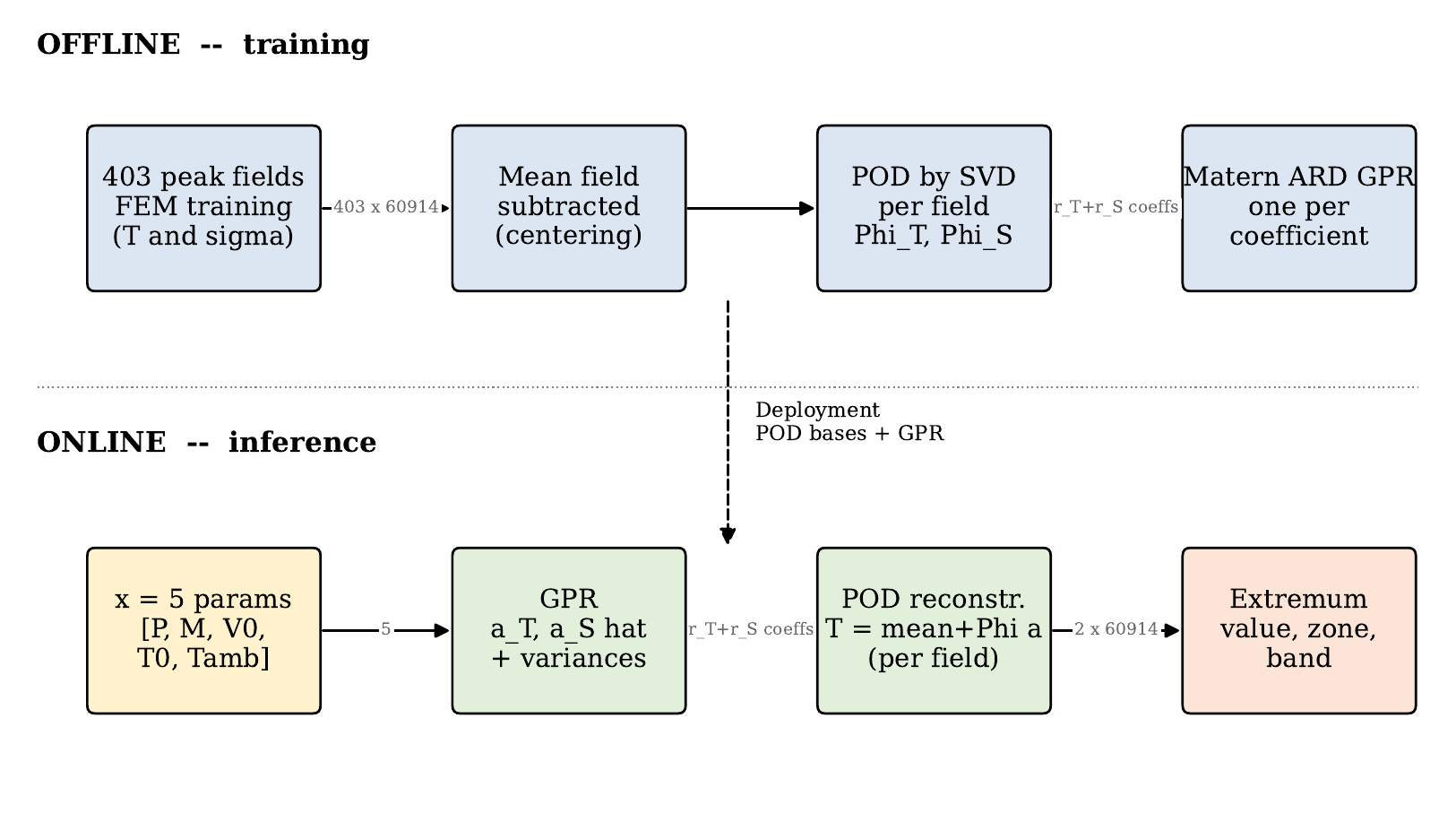}
  \caption{The non-intrusive POD-GPR pipeline. During the offline stage, we construct a proper orthogonal decomposition basis for each field and train a Gaussian process per reduced coordinate. During the online stage, the processes predict coordinates for new operating conditions, allowing full-field reconstruction and peak evaluation.}
  \label{fig:pipeline}
\end{figure}

The emulator maps parameters directly to the peak fields, bypassing solver equations and operating strictly as a data-driven model. This non-intrusive approach requires only input-output snapshots rather than full access to governing differential operators, contrasting with classical projection-based reduced-order schemes. Figure~\ref{fig:pipeline} outlines the two distinct phases of the formulation. During the offline phase, we construct a reduced-order basis for each physical field and train a regression model to map inputs to coordinates in this lower-dimensional space. During the online phase, we query the trained regressor at any new parameter vector $\mathbf{x}$ to reconstruct both physical fields. We perform the spatial reduction via proper orthogonal decomposition (POD) and resolve the regression step using a Gaussian process in the coefficient space to obtain simultaneous predictions and associated variances. Although both reductions run independently, they share the same spatial discretization over the fixed finite-element mesh.

\subsection{Proper orthogonal decomposition}
\label{subsec:pod}

We reduce each physical field through proper orthogonal decomposition. This process begins by centering the $\nTrain{}$ training snapshots about their mean field and organizing them into a snapshot matrix. Singular value decomposition of this matrix yields spatial modes and associated singular values, where the squared singular values represent the modal energy~\cite{sirovich1987snapshots,berkooz1993pod}. These modes define the directions of highest spatial variation across the training data~\cite{chatterjee2000podintro}. We define the inner product using the Euclidean product of the nodal vectors, which fits our near-uniform mesh layout. Strongly graded meshes would instead require a mass-matrix-weighted inner product to account for variable element volumes. Utilizing the method of snapshots allows us to solve the underlying eigenproblem at a computational cost scaled to the training set size instead of the mesh size~\cite{sirovich1987snapshots}. The resulting projection onto the first $r$ modes yields the low-dimensional field approximation:
\begin{equation}
  u \approx \bar{u} + \Phi_r\, a,
  \qquad
  a = \Phi_r^{\top}(u - \bar{u}),
  \label{eq:pod-projection}
\end{equation}
where $\bar{u}$ is the mean field, $\Phi_r$ contains the first $r$ spatial modes, and $a \in \mathbb{R}^{r}$ represents the vector of reduced coordinates~\cite{rathinam2003podanalysis}.

\subsection{Rank selection by peak accuracy}
\label{subsec:rank}

Truncating the proper orthogonal decomposition basis minimizes the spatial energy norm of the reconstruction error rather than the localized peak error. These two error metrics behave differently in this setting. A low global energy error does not guarantee accurate peak predictions. The peak is a localized, low-energy feature that leading modes often fail to resolve. For instance, the first mode accounts for \modeUnPicS{} of the total energy for the von Mises field and \modeUnPicT{} for the thermal field. While the cumulative energy variance saturates using only \rankPicNNNS{} modes for the von Mises field, the localized peak reconstruction error continues to decline at higher ranks. Selecting the truncation limit based strictly on cumulative energy variance can lead to an under-resolved peak. We address this by determining the truncation rank directly from the validation set peak reconstruction error. We select the minimum rank where this peak error falls safely below our target threshold, yielding a working rank of $\rangOpere{}$ for both fields. This empirical selection remains subject to the formal verification and exact error bounds derived in Section~\ref{sec:budget}.

\subsection{Gaussian process regression of the reduced coordinates}
\label{subsec:gpr}

We map the parametric inputs to the low-dimensional coordinates using Gaussian process regression~\cite{rasmussen2006gpml}. We train an independent Gaussian process for each coordinate, requiring $\rangOpere{}$ processes per field and $2\,\rangOpere{}$ processes for the two fields in total. Treating these coordinates as independent neglects cross-modal correlations and the physical coupling between thermal and mechanical fields. This standard decoupling simplification allows efficient training, while the resulting modeling errors are captured directly by the regression term in the error budget of Section~\ref{sec:budget}.

Each process maps the input vector $\mathbf{x} \in \mathbb{R}^{\nParam}$ to a single standardized coordinate to yield a posterior mean and a corresponding predictive variance. The covariance structure uses a Mat\'ern kernel~\cite{stein1999kriging} with automatic relevance determination, assigning an independent length scale to each input parameter. We optimize these hyper-parameters by maximizing the log-marginal likelihood using multiple random restarts. Adding a minor numerical jitter to the diagonal of the covariance matrix stabilizes the linear algebra, rendering the processes near-interpolating with minimal noise. The mean field $\bar{u}$ from \eqref{eq:pod-projection} acts as a static offset, meaning each Gaussian process operates with a constant prior mean in the standardized space. The optimized length scales automatically scale down input dimensions that have minor influence on the coordinate predictions.

During online evaluation at a new parameter point $\mathbf{x}$, each independent process predicts its corresponding coordinate, we reconstruct the full field via \eqref{eq:pod-projection}, and we apply the spatial maximum operator defined in \eqref{eq:peak-functional}. This architecture bypasses assembling or solving the finite-element system, qualifying the framework as a fully non-intrusive surrogate model~\cite{guo2018podgpr,ortali2022podgprfluid}. Alternative regression techniques include feedforward neural networks~\cite{hesthaven2018neuralnetworks} or manifold interpolation methods~\cite{amsallem2008interpolation}. In their standard deterministic form these do not natively output the calibrated predictive variance necessary for the a posteriori peak error budget in Section~\ref{sec:budget} and the level-one uncertainty of Section~\ref{sec:results}, although they can be equipped with uncertainty estimates by other means.
\section{Peak error budget and level-one uncertainty}
\label{sec:budget}

We evaluate the non-intrusive emulator against the high-fidelity finite-element reference at the peak state. This assessment uses an a posteriori error budget bounded by the numerical discretization floor established in Section~\ref{subsec:floor}. We also compute a level-one predictive uncertainty for the certified peak. This evaluation operates on a surrogate-to-code basis, demonstrating that the actual emulator peak error falls within the discretization uncertainty of the reference solver.

\subsection{A two-term peak error budget}
\label{subsec:two-term}

The certified quantity of interest is the spatial peak of each field, $Q(u) = \max_{1 \le i \le N_{\mathrm{n}}} u_i$, defined in Section~\ref{sec:problem}. Because the peak operator is nonlinear, the spatial node where this maximum occurs can vary among the high-fidelity, projected, and emulated fields. We decompose the total prediction error on this scalar quantity into two distinct sources: proper orthogonal decomposition truncation and reduced coordinate regression. Let $Q_{\mathrm{ref}} = Q(u_{\mathrm{ref}})$ denote the peak value of the reference solver field. We define the projected peak as $Q_{\mathrm{proj}} = Q\bigl(\bar{u} + \Phi_r \Phi_r^{\top}(u_{\mathrm{ref}} - \bar{u})\bigr)$, which represents the spatial maximum obtained by reconstructing the field from its exact projected coordinates. Let $Q_{\mathrm{emu}} = Q(u_{\mathrm{emu}})$ be the peak value predicted by the emulator. The absolute realized error, the basis truncation error, and the coordinate regression error are formulated as:
\begin{equation}
  e_{\mathrm{real}} = \bigl\lvert Q_{\mathrm{emu}} - Q_{\mathrm{ref}} \bigr\rvert,
  \quad
  e_{\mathrm{trunc}} = \bigl\lvert Q_{\mathrm{proj}} - Q_{\mathrm{ref}} \bigr\rvert,
  \quad
  e_{\mathrm{reg}} = \bigl\lvert Q_{\mathrm{emu}} - Q_{\mathrm{proj}} \bigr\rvert.
  \label{eq:budget-terms}
\end{equation}
The truncation term isolates the geometric error of the reduced basis without incorporating regression inaccuracies. The regression term captures the error introduced solely by the Gaussian process models on this fixed basis. We compute these metrics using the independent test set containing $\nTest{}$ cases.

\subsection{Error attribution}
\label{subsec:attribution}

The realized prediction error on the scalar peak represents an exact algebraic sum of the signed truncation and regression contributions:
\begin{equation}
  Q_{\mathrm{emu}} - Q_{\mathrm{ref}}
    = \bigl(Q_{\mathrm{proj}} - Q_{\mathrm{ref}}\bigr)
    + \bigl(Q_{\mathrm{emu}} - Q_{\mathrm{proj}}\bigr).
  \label{eq:signed-decomposition}
\end{equation}
Applying the triangle inequality yields a rigorous upper bound on the absolute realized peak error: $e_{\mathrm{real}} \le e_{\mathrm{trunc}} + e_{\mathrm{reg}}$. Case-by-case evaluation of these absolute terms across the test set confirms that this bound holds without exception. This formulation attributes the prediction error directly to its underlying mathematical sources. The basis truncation error is the primary driver of the total error, dominating the regression term in \budDomFracT{} of the thermal test cases and \budDomFracS{} of the von Mises stress cases. The linear dimensional reduction governs the prediction accuracy, while the Gaussian process approximation contributes minimally. Section~\ref{sec:results} reports the relative peak errors.

\subsection{Bounding against the discretization floor}
\label{subsec:below-floor}

Evaluating the realized peak error requires comparison with the baseline discretization floor of the solver established in Section~\ref{subsec:floor}~\cite{roy2011vvuq}. A point-by-point comparison shows that the actual emulator peak error falls below this discretization floor in \budSousPlancherFrac{} of the test cases for both physical fields. The surrogate reproduces the high-fidelity solver within its own numerical resolution limits. This alignment defines and limits our surrogate-to-code verification because attempting to achieve errors below this mesh-dependent threshold lacks physical justification.

\subsection{Level-one uncertainty of the peak}
\label{subsec:uq}

We propagate the posterior coordinate variances computed by the Gaussian process regression models directly to the peak quantities of interest~\cite{rasmussen2006gpml}. Using a Monte Carlo sampling procedure, we draw coordinate realizations from the posterior distributions, reconstruct the spatial fields, and extract the resulting maxima. The empirical quantiles of this generated peak distribution define a level-one predictive interval, isolating the uncertainty of the scalar peak from the full-field spatial uncertainty. This interval measures regression uncertainty alone, excluding basis truncation errors and solver discretization limits, which are handled independently within our error budget.

This propagation method exhibits two distinct mathematical behaviors. First, our Monte Carlo framework treats the coordinate posteriors as independent. The maximum off-diagonal correlation of the validation residuals reaches only \corrHorsDiagS{} for the von Mises coordinates, and incorporating a full-covariance structure produces negligible changes in the peak coverage, verifying the validity of the diagonal assumption. Second, the maximum of a finite sample represents a downward-biased estimator of the expected peak. The predictive interval remains slightly conservative, acting as a practical engineering bound rather than a strictly unbiased statistical estimator.

We evaluate the calibration of this predictive interval on the validation split at a nominal level of \nomLevel{}. The actual spatial peak falls within the predicted interval in \covExtValT{} of the thermal cases and \covExtValS{} of the von Mises stress cases. The thermal predictions display slight over-coverage, while the mechanical stress predictions exhibit a minor coverage deficit. This deficit remains statistically insignificant across the validation sample because the corresponding Clopper-Pearson confidence interval, [\cpLoExtValS{},\, \cpHiExtValS{}], encompasses the nominal target. The computed level-one uncertainty remains consistent with the nominal specifications, avoiding claims of exact calibration.
\section{Results on the held-out test set}
\label{sec:results}

\subsection{Realized peak accuracy at the quantity of interest}
\label{subsec:accuracy}
\begin{figure}[htbp]
  \centering
  \includegraphics[width=\linewidth]{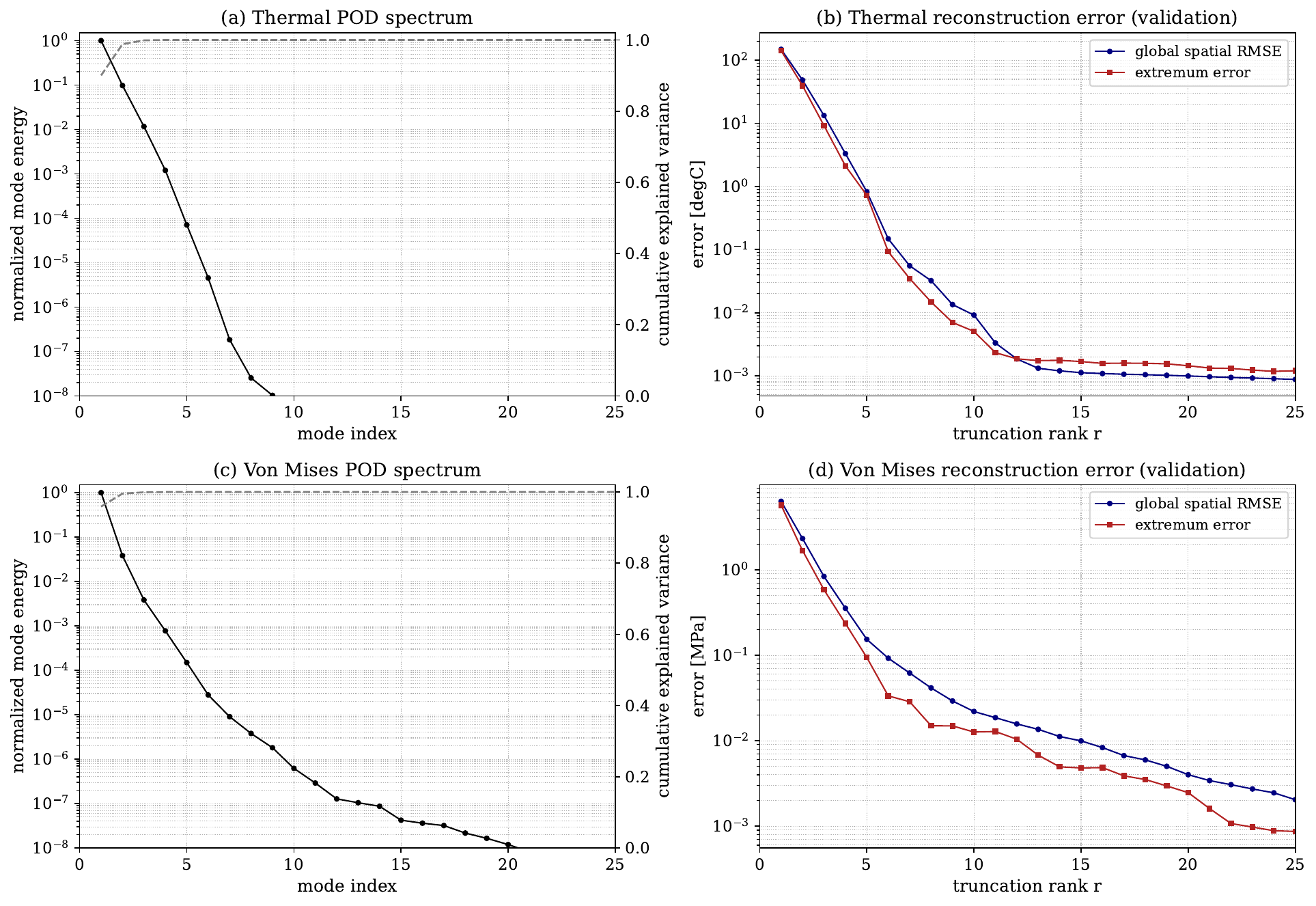}
  \caption{Reconstruction error of each field against the retained rank. The peak reconstruction error keeps decreasing beyond the rank at which the cumulative variance saturates, which is why the rank is selected on the peak error rather than on the variance.}
  \label{fig:pod-spectrum}
\end{figure}

\begin{figure}[htbp]
  \centering
  \includegraphics[width=\linewidth]{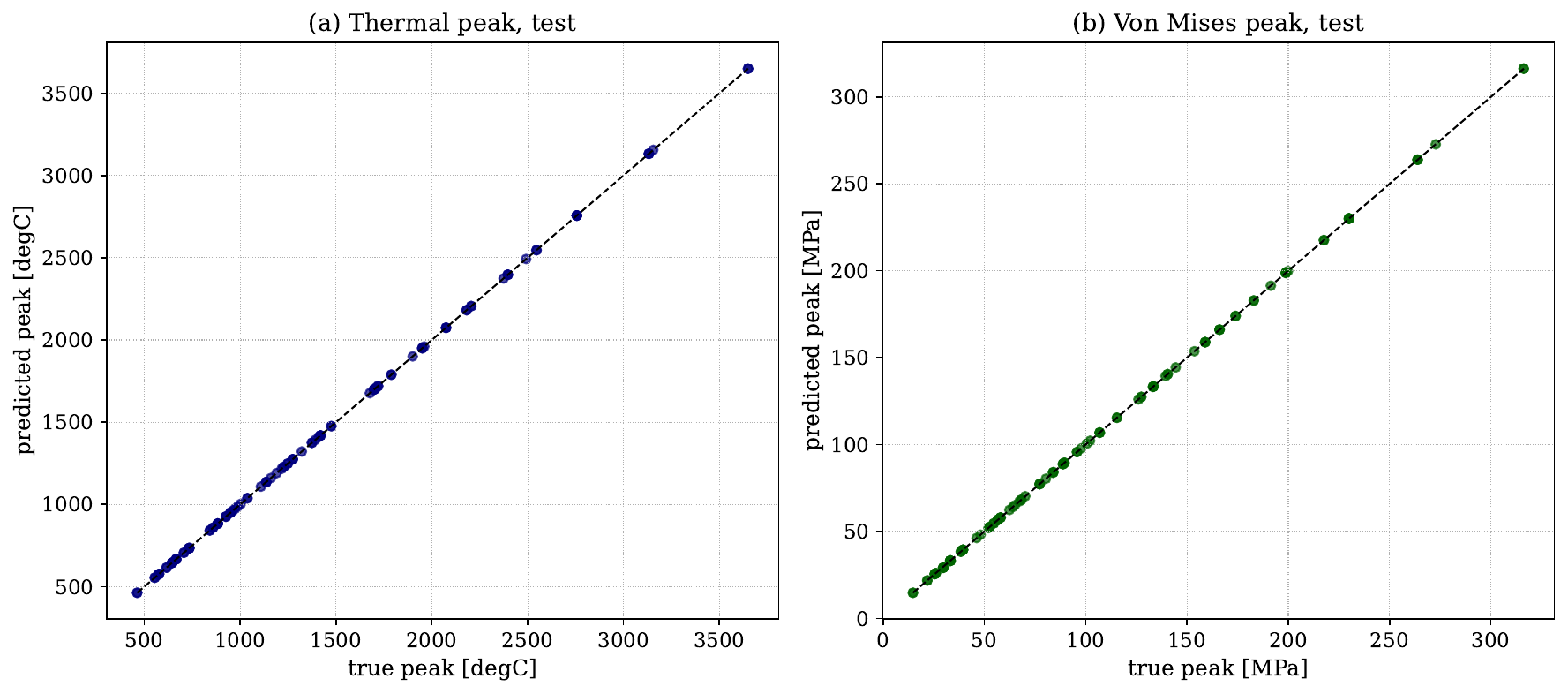}
  \caption{Predicted against reference peak values on the held-out test set, for the thermal peak and the von Mises peak. Points on the diagonal indicate exact agreement, and the band corresponds to the pre-registered conformity criterion.}
  \label{fig:parity}
\end{figure}

\begin{figure}[htbp]
  \centering
  \includegraphics[width=\linewidth]{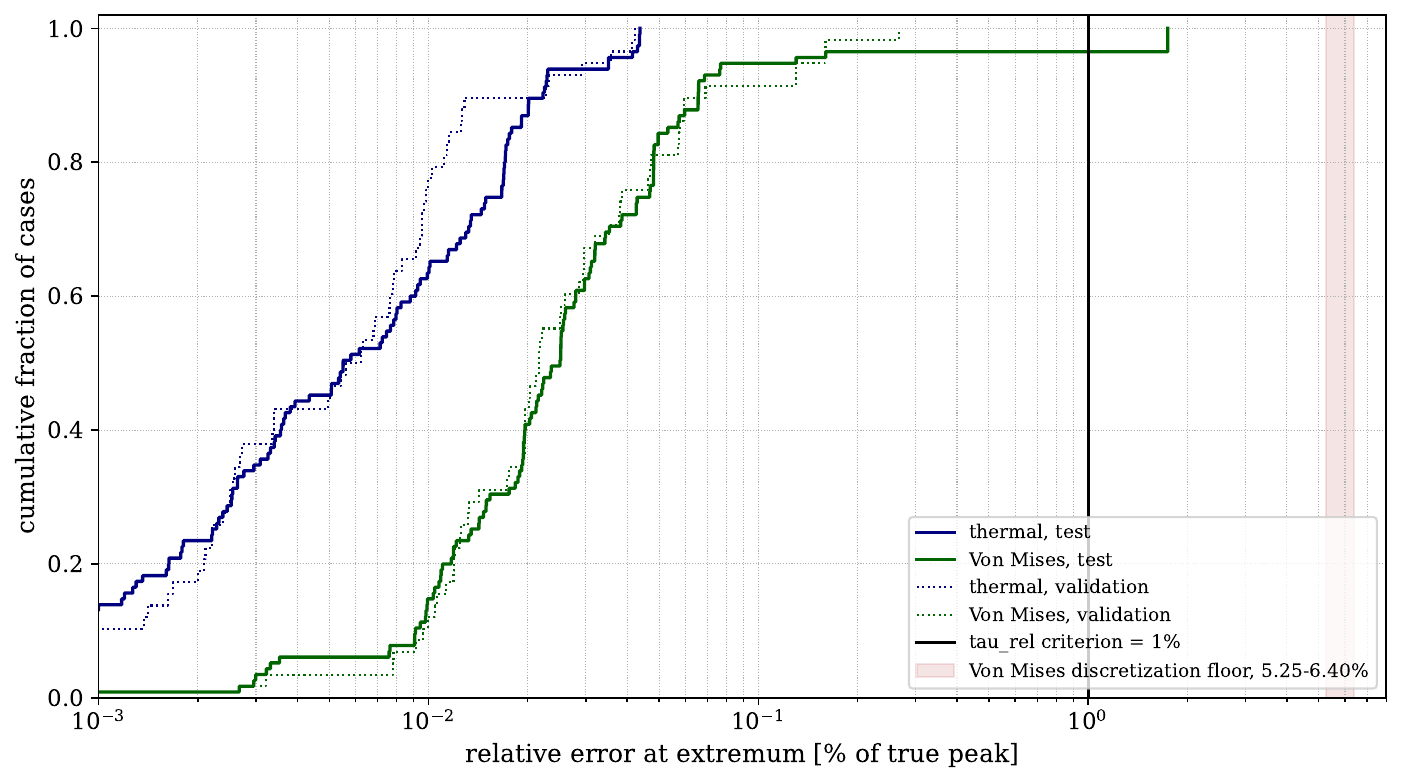}
  \caption{Cumulative distribution of the relative peak error on the held-out test set, for both fields, on a logarithmic axis. The discretization floor is marked, and the bulk of the distribution lies below it.}
  \label{fig:error-cdf}
\end{figure}

We evaluate the emulator's predictive performance using the spatial extremes of the temperature and von Mises stress fields over the production mesh consisting of $N_{\mathrm{n}} = \nNodes{}$ nodes. To ensure a rigorous assessment, we report all subsequent metrics on an independent, held-out test set of $\nTest{}$ cases. This test partition remained strictly blinded until we finalized the rank-$\rangOpere{}$ emulator design, preventing any post-hoc tuning or optimization. We define a compliance threshold using a relative peak error limit of $\tauRel{}$, which corresponds to the structural design tolerances of the physical brake assembly.

For the thermal peak, the emulator achieves a median relative error of $\medRelExtTestT{}$ and a maximum relative error of $\maxRelExtTestT{}$ (corresponding to a maximum absolute deviation of $\absExtTestT{}$), yielding a thermal compliance rate of $\conformFracTestT{}$. For the peak von Mises stress, the median relative error is $\medRelExtTestS{}$ and the maximum is $\maxRelExtTestS{}$ (representing an absolute deviation of $\absExtTestS{}$). This results in a mechanical compliance rate of $\conformFracTestS{}$, where only $\nHorsCritTestS{}$ out of $\nTest{}$ test cases exceed the target limit. We report these unadjusted compliance rates in strict accordance with our blind-testing protocol. Because the full stress field is not fully converged at a truncation rank of $\rangOpere{}$, peak identification on the truncated basis carries an epistemic risk under parametric variation. Figure~\ref{fig:parity} displays the parity plots of predicted versus reference peak values, and Figure~\ref{fig:error-cdf} shows the cumulative distribution of the relative errors.

These performance metrics depend on two physical factors. First, the relative tolerance limit of $\tauRel{}$ is highly stringent, yet the few stress cases exceeding this threshold remain significantly below the numerical solver's discretization error described in Section~\ref{subsec:floor-results}. Second, because we selected the truncation rank to minimize peak error rather than global reconstruction energy, peak accuracy converges at a lower modal dimension than the full-field solution. Consequently, while the global von Mises field reconstruction does not reach a complete plateau at rank $\rangOpere{}$ (Figure~\ref{fig:pod-spectrum}), the emulator maintains high fidelity at the certified spatial peaks. We address the spatial localization of these extreme values in Section~\ref{sec:ablation}.

\subsection{The two-term peak error budget}
\label{subsec:budget-results}
\begin{figure}[!htbp]
  \centering
  \includegraphics[width=\linewidth]{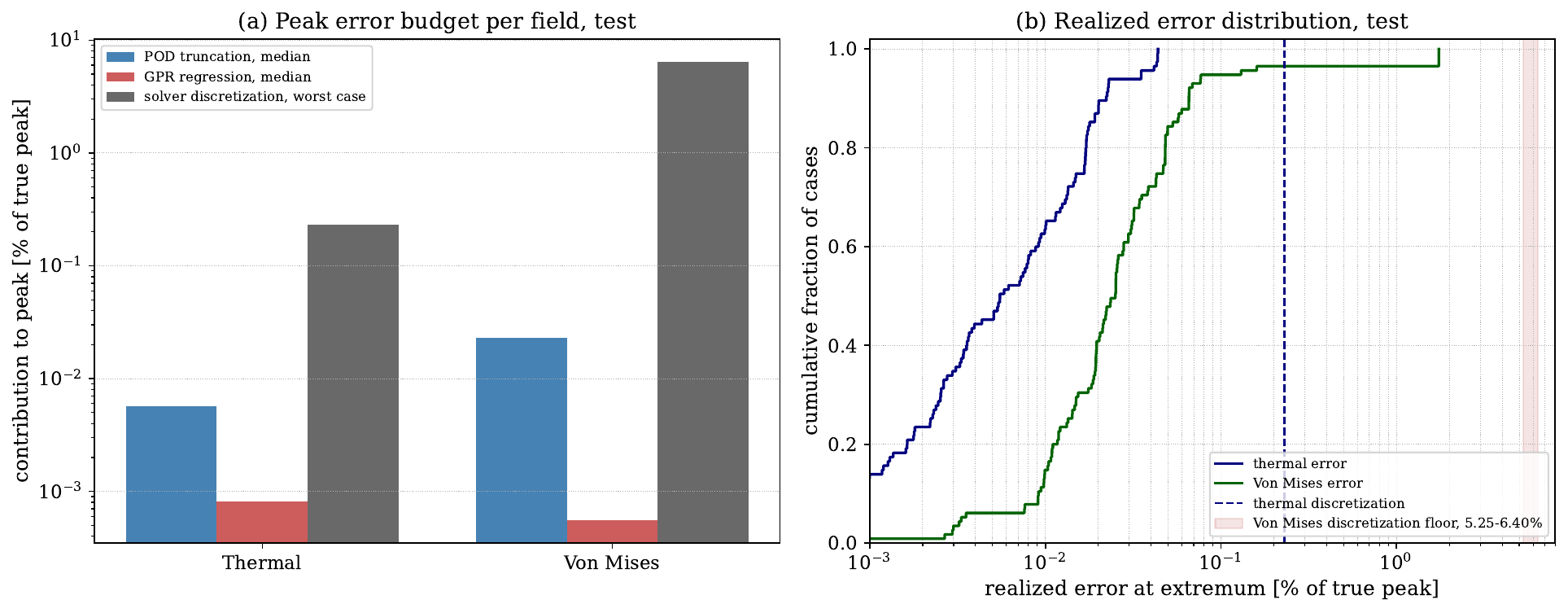}
  \caption{Two-term peak error budget on the held-out test set. For each field the realized peak error is shown against the truncation term and the regression term, with the discretization floor as a shaded band. The realized error lies below the floor, and the truncation term dominates the regression term.}
  \label{fig:budget}
\end{figure}

The total realized error on the scalar peak for any given test configuration equals the algebraic sum of the signed truncation and regression errors. The triangle inequality guarantees that the sum of these absolute components bounds the total absolute error, a condition verified across all $\nTest{}$ test cases. The finite-element discretization error is not a term of this budget. It is reported separately as a solver reference floor in Section~\ref{subsec:floor-results}.

In physical units, the thermal truncation error exhibits a median of $\budTruncMedT{}$ and a $95\textsuperscript{th}$ percentile of $\budTruncPToT{}$, compared to a regression error median of $\budRegMedT{}$ and a $95\textsuperscript{th}$ percentile of $\budRegPToT{}$. The total realized thermal error has a median of $\budRealMedT{}$ and a $95\textsuperscript{th}$ percentile of $\budRealPToT{}$. For the peak von Mises stress, the truncation error has a median of $\budTruncMedS{}$ and a $95\textsuperscript{th}$ percentile of $\budTruncPToS{}$, while the regression error has a median of $\budRegMedS{}$ and a $95\textsuperscript{th}$ percentile of $\budRegPToS{}$, producing a total realized stress error median of $\budRealMedS{}$ and a $95\textsuperscript{th}$ percentile of $\budRealPToS{}$. Truncation errors exceed regression errors in $\budDomFracT{}$ of the thermal simulations and $\budDomFracS{}$ of the mechanical stress simulations. The spatial reduction step governs the overall emulator accuracy, suggesting that performance improvements depend on expanding the linear basis rather than optimizing the regression hyperparameters. Figure~\ref{fig:budget} displays this case-by-case error distribution.

\subsection{The discretization floor and the surrogate-to-code fidelity}
\label{subsec:floor-results}
\begin{figure}[htbp]
  \centering
  \includegraphics[width=\linewidth]{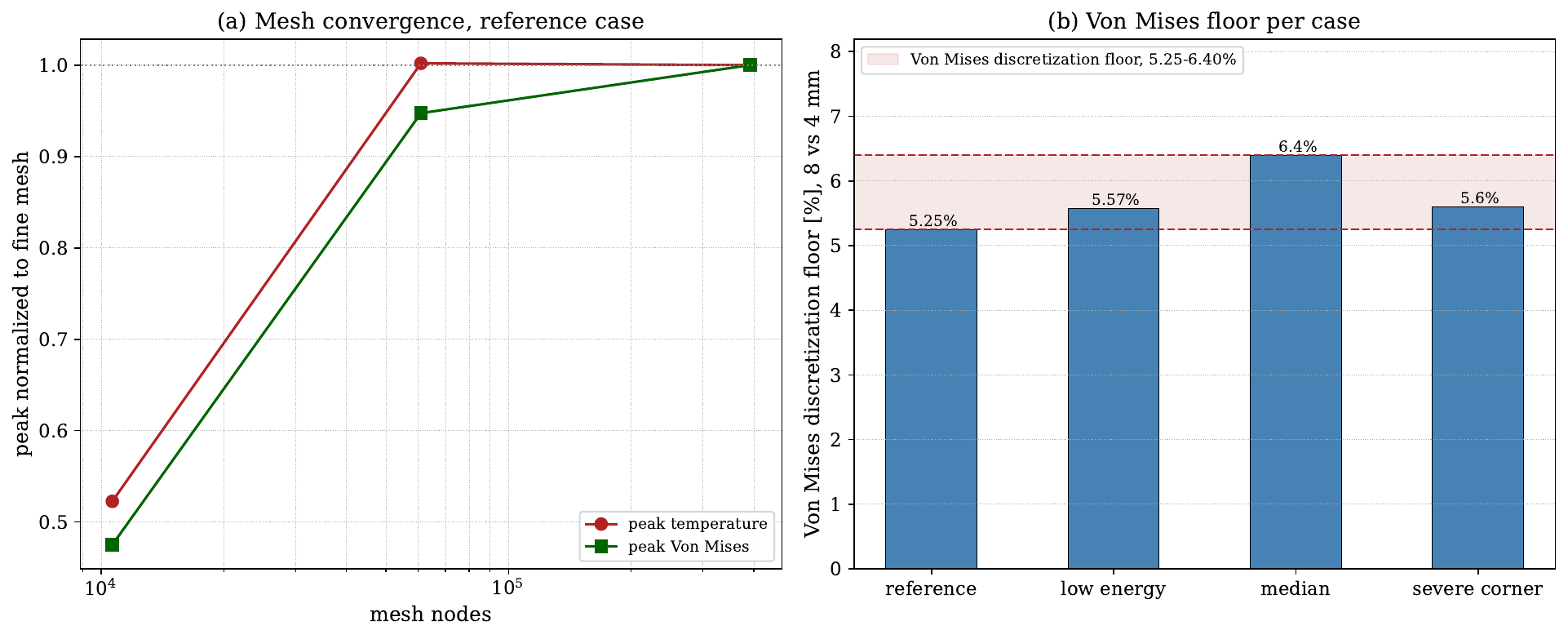}
  \caption{Three-mesh convergence of the peak temperature and the peak von Mises stress. The peak values are shown against the mesh size, from the coarse mesh to the reference mesh. The thermal peak is converged by the production mesh, while the von Mises peak retains a mesh dependence that defines the per-case discretization floor.}
  \label{fig:mesh-convergence}
\end{figure}

We evaluate the emulator accuracy against the finite-element discretization floor defined in Section~\ref{subsec:floor} and illustrated via the mesh convergence study in Figure~\ref{fig:mesh-convergence}. This analysis compares a coarse mesh of $\meshCoarseSize{}$ mm ($\meshCoarseNodes{}$ nodes), which under-resolves the stress field with a peak von Mises value of $\svmMaxCoarse{}$ MPa, against our production mesh of $\meshProdSize{}$ mm ($\meshProdNodes{}$ nodes) and a highly refined reference mesh of $\meshFineSize{}$ mm ($\meshFineNodes{}$ nodes). While the thermal peaks converge fully on the production mesh, the peak von Mises stress exhibits residual mesh sensitivity. We evaluate this discretization floor, defined as the relative difference between the production and refined meshes, on a case-by-case basis. The stress discretization floor spans from $\epshRelSMin{}$ to $\epshRelSMax{}$. The maximum relative discretization error occurs near the median loading cases rather than the extreme loading cases because the mesh sensitivity correlates with the spatial gradient and location of the peak stress rather than its absolute magnitude. In contrast, the thermal discretization floor remains negligible, ranging from $\epshRelTMedian{}$ to $\epshRelTCoin{}$.

The total realized emulator error remains below the numerical discretization floor in $\budSousPlancherFrac{}$ of the test cases across both physical fields. Specifically, the worst-case relative stress error of $\maxRelExtTestS{}$ stays below the minimum stress discretization floor of $\epshRelSMin{}$. When the surrogate model's approximation error drops below the underlying discretization uncertainty of the solver, further optimization of the emulator fails to yield physically or numerically meaningful improvements. This is the basis for reading the fidelity as surrogate-to-code, under the present constant-property, radiation-free reference model and fixed mesh family. This assessment measures the emulator's ability to replicate the numerical solver rather than the physical brake assembly, which the underlying finite-element formulation tends to overestimate.

\subsection{Computational speedup}
\label{subsec:speedup}

The surrogate framework replaces the high-fidelity finite-element solver (mean runtime of $\femMoyS{}$ seconds) with a rapid evaluation step (median runtime of $\tInfMedMs{}$ ms, excluding disk I/O for loading bases and Gaussian process hyperparameters). Since the numerical solver's execution time scales with the physical braking duration, the computational speedup spans a wide operational range. Using a single reference machine to ensure consistent timing benchmarks, the local speedup factor ranges from $\accelMin{}$ to $\accelMax{}$, with a median acceleration of $\accelMed{}$. Over this benchmark suite, the total solver time of $\femCumulS{}$ seconds is reduced by a factor of $\accelCumul{}$. We provide the detailed distribution and range of these speedup factors in Figures~\ref{fig:speedup-distribution} and \ref{fig:speedup-range} in \ref{app:supplementary}.

\subsection{Off-grid robustness and the bound on validity}
\label{subsec:offgrid}

We evaluate the robustness of the frozen emulator on $\genNHgPoints{}$ Latin-hypercube samples generated within the parameter boundaries but off the structured factorial grid. For these off-grid test points, the thermal peak prediction error exhibits a median of $\genMedRelHgT{}$ and a maximum of $\genMaxRelHgT{}$, maintaining a compliance rate of $\genConfHgT{}$. The peak von Mises stress error shows a median of $\genMedRelHgS{}$ and a maximum of $\genMaxRelHgS{}$, corresponding to a compliance rate of $\genConfHgS{}$. Compared to the on-grid test partition, the median thermal peak error increases from $\medRelExtTestT{}$ to $\genMedRelHgT{}$ and the compliance rate decreases from $\conformFracTestT{}$ to $\genConfHgT{}$. This moderate degradation reflects the challenges of off-grid extrapolation within the bounded input domain. This performance gap is more pronounced in the temperature field than in the stress field, an observation whose mechanism we do not establish here and take up in Section~\ref{sec:discussion}.

These metrics apply strictly within the boundaries of the sampled parameter space. This off-grid analysis serves as a secondary check to test model robustness and verify the operational limits defined in Section~\ref{sec:hifi}. The size of the off-grid dataset ($\genNHgPoints{}$ points) remains constrained by the computational cost of the high-fidelity solver rather than the surrogate, making these results indicative of model behavior under off-grid conditions.

\subsection{Summary of the test-set assessment}
\label{subsec:results-summary}

The independent test-set evaluation yields three key conclusions. First, the emulator meets the target structural compliance criteria at the spatial peaks for the thermal field and the vast majority of the mechanical stress configurations. Second, basis truncation represents the primary contributor to the total peak error, whereas regression errors remain secondary. Third, the total prediction error remains below the solver's inherent discretization floor across all test cases, justifying our surrogate-to-code verification framework.
\section{Ablation, scope of the linear reduction, and interpretability}
\label{sec:ablation}

We analyze the primary sources of emulator approximation error, evaluate the physical limits of linear dimensional reduction, and interpret the physical dependencies captured by the surrogate model.

\subsection{Ablation study and error sources}
\label{subsec:ablation}

The two-term error budget in Section~\ref{sec:budget} attributes the peak prediction error to the proper orthogonal decomposition basis rather than the Gaussian process regression. Our ablation analysis confirms this attribution, particularly in the worst-performing test cases where truncation error dominates. For the minority of test cases where the reference peak von Mises stress falls outside the level-one predictive interval from Section~\ref{subsec:uq}, the regression error remains negligible compared to the truncation error. These specific cases correspond to high braking loads coupled with low aircraft masses, a regime where the localized spatial gradients are under-resolved at the selected basis rank $r$. Consequently, the linear dimensional reduction step acts as the primary accuracy bottleneck rather than the Gaussian process regression.

\subsection{Scope of the linear reduction}
\label{subsec:scope}

We assess whether adopting a nonlinear dimensional reduction scheme would improve peak prediction accuracy. Within the analyzed parameter bounds, a nonlinear reduction does not yield measurable improvements because the POD truncation error already falls below the numerical discretization floor established in Section~\ref{subsec:floor}. Any further reduction in reconstruction error would be obscured by the underlying finite-element mesh resolution, rendering more complex reduction techniques unobservable in this verification framework. The theory of the Kolmogorov $n$-width offers a suggestive but not conclusive backdrop~\cite{greif2019nwidth}, since it applies rigorously to linear coercive problems, whereas the present problem couples thermal diffusion with a thermoelastic response. The linear $n$-width of the purely diffusive thermal field decays rapidly, though we do not assume the coupled thermoelastic equations satisfy strict coercivity. Nonlinear reduction methods (including kernel principal component analysis~\cite{scholkopf1998kpca}, autoencoders~\cite{lee2020autoencoder,fresca2021dlrom}, and their kernel-mapped or Gaussian-process extensions~\cite{salvador2021kpodnn,deshpande2025gpae}) are designed for transport-dominated physical processes where the linear $n$-width decays slowly. Because our transient thermomechanical problem is not transport-dominated, we state the superfluity of a nonlinear reduction as a scoped corollary, valid for this problem class and parameter domain, and not as a general claim.

\subsection{Zonal localization of the extremum}
\label{subsec:zone}

The secondary contribution predicts the extremum as a zone, the super-level set $Z_\alpha$ of Section~\ref{sec:problem} at an operating threshold $\alpha$, rather than as a node. We state its scope with care, since a trivial zone would carry no information. For the thermal field, this super-level set is spatially trivial at the selected threshold, spanning a complete annulus over the friction band with negligible variation across parameter cases. For the von Mises stress field, the spatial centroid and radial bounds of the zone vary less than the characteristic zone thickness, meaning the physical location of the peak zone remains too insensitive to serve as a reliable predictive output. In contrast, the total surface area of the high-stress zone represents a stable, physically informative metric. The mean area of this stress zone is $\zoneAireMoyS{}$ of the total friction band surface, exhibiting a case-to-case coefficient of variation of $\zoneAireCVS{}$. This surface area displays strong monotonic correlation with aircraft mass, braking load, and initial speed across the training and validation subsets. The emulator successfully captures this physical dependency, achieving a correlation coefficient of $\zoneRhoAireS{}$ on the validation set. The initial and ambient temperatures exert negligible influence on the spatial zone size. The weak correlation between zone area and initial temperature observed in the validation subset represents a finite-sample artifact present in the underlying finite-element solver data. That the emulator reproduces this local trend reflects surrogate-to-code tracking of the observed data rather than a physical causal relationship, and it remains bounded by the discretization floor, so it does not degrade the surrogate-to-code fidelity. This spatial area provides a stable scalar indicator of stress concentration that avoids the numerical instabilities associated with tracking a single maximum node.

\subsection{Interpretability through the length scales}
\label{subsec:interpret}
\begin{figure}[htbp]
  \centering
  \includegraphics[width=\linewidth]{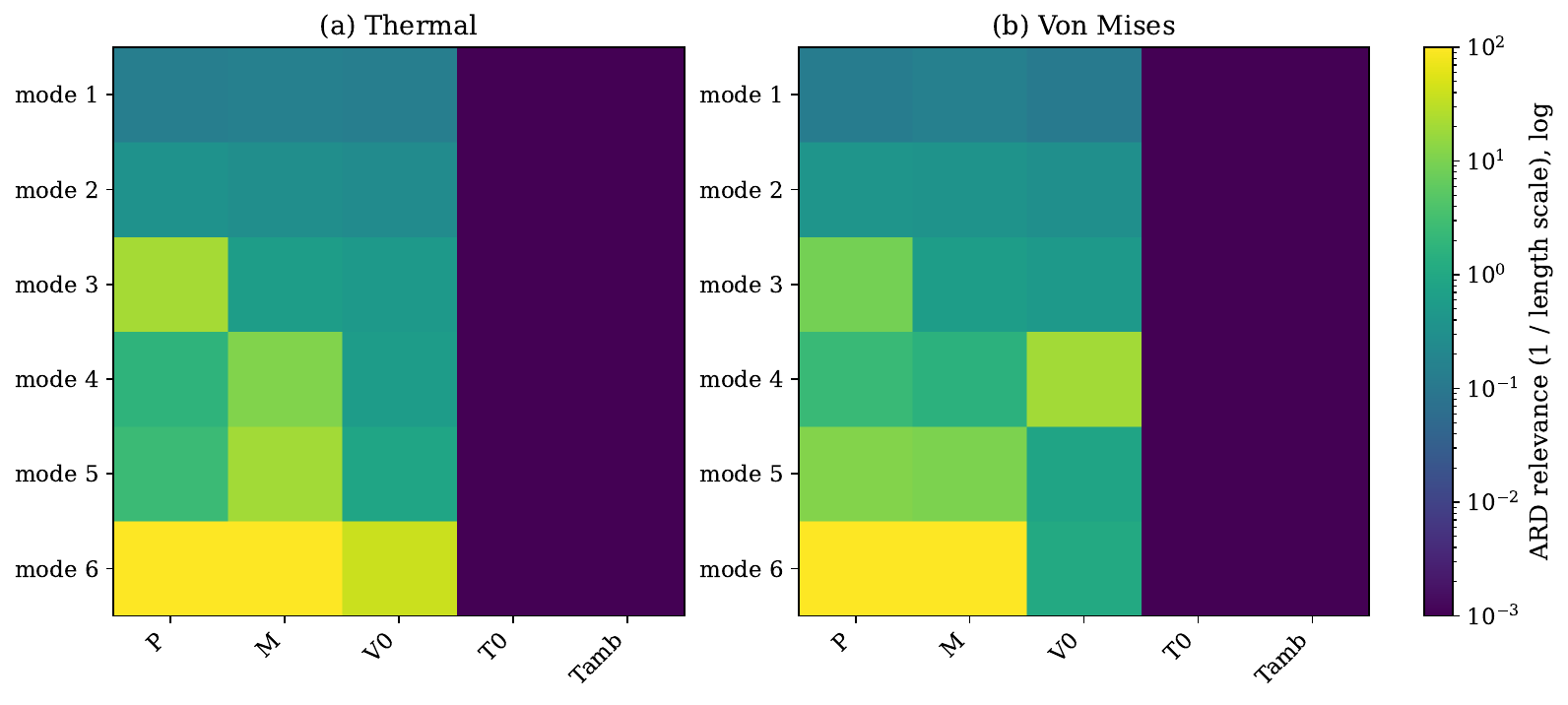}
  \caption{Fitted automatic-relevance-determination length scales, one per operating parameter. A long length scale marks a parameter to which the reduced coordinates vary slowly, and thus a low influence. The initial and ambient temperatures are the least influential, the braking load, the mass, and the initial speed the most.}
  \label{fig:ard-relevance}
\end{figure}

The automatic relevance determination (ARD) framework of the Gaussian process regression fits an independent correlation length scale $\ell_d$ to each input dimension. Under the chosen Mat\'ern kernel, these length scales provide a qualitative measure of parameter influence. Larger length scales indicate that the predicted coordinates vary slowly along the corresponding parameter dimension, signifying a lower impact on the output. Based on this analysis, the initial and ambient temperatures represent the least influential parameters, matching their minimal impact on the high-stress zone area. Conversely, the braking load, aircraft mass, and initial velocity act as the primary drivers of the physical response, which is consistent with the energy dissipation and mechanical loading mechanisms detailed in Section~\ref{sec:hifi}. Figure~\ref{fig:ard-relevance} displays these fitted length scales. These ARD parameters are model-dependent and differ from variance-based Sobol sensitivity indices, though both metrics yield comparable qualitative rankings~\cite{saltelli2008primer,sobol2001indices}. While Gaussian process models can be post-processed to calculate formal Sobol indices~\cite{marrel2009sobolgp}, and polynomial chaos expansions natively support global sensitivity analysis~\cite{xiu2002wienaskey,sudret2008pcegsa}, we present the ARD length scales strictly as a qualitative interpretation of parameter importance.

\subsection{Section summary}
\label{subsec:ablation-summary}

Our ablation analysis establishes that proper orthogonal decomposition truncation is the dominant source of approximation error in the emulator pipeline. We verify that a linear reduced-order basis remains sufficient within the constraints of the solver's discretization floor, and we demonstrate that the non-intrusive surrogate accurately reproduces the physical sensitivities governing the spatial extent of the stress concentrations.
\section{Discussion and limitations}
\label{sec:discussion}

This section contextualizes the verification results, delineates the physical and numerical boundaries of the emulator, and positions the framework relative to existing reduced-order modeling methodologies.

\subsection{What the assessment establishes}
\label{subsec:establishes}

Our evaluation establishes a verified error budget anchored directly to the solver discretization floor within a formal verification and validation framework~\cite{roy2011vvuq}. Across the independent test set, the emulator replicates the high-fidelity solver at the certified spatial peaks within the numerical discretization uncertainty of the reference model. This surrogate-to-code agreement verifies the emulator against the underlying numerical solver rather than validating it against physical experimental measurements. 

Anchoring the error budget to this numerical floor establishes a rigorous stopping criterion, as optimizing the emulator beyond the spatial resolution of the reference solver yields no physical or mathematical utility. Once the emulator's approximation error falls below this discretization uncertainty, further parameter tuning offers no actionable improvement in the fidelity of that numerical model. The two-term budget attributes this residual error. By evaluating the regression error on the fixed, truncated basis, we demonstrate that spatial truncation, rather than Gaussian process regression, governs the baseline surrogate error.

\subsection{Limitations}
\label{subsec:limitations}

We identify several key limitations that define the scope of these findings.

First, the peak von Mises stress reference is not fully converged on the production mesh. The stress discretization floor consequently acts as a case-dependent range rather than a static scalar, and we assert emulator compliance relative to this local interval. Although the realized peak errors on the test set consistently fall below the conservative lower bound of this range, establishing a fully mesh-independent stress reference would require resolving all training snapshots on the refined mesh, a step omitted due to computational budget constraints.

Second, the level-one uncertainty derived in Section~\ref{subsec:uq} accounts strictly for epistemic regression uncertainty. It excludes the independent truncation and discretization errors, meaning the resulting predictive interval does not represent the full statistical uncertainty of the physical peak. The empirical stress coverage of $\covExtValS{}$ slightly underperforms the nominal $\nomLevel{}$ target, a discrepancy that remains statistically insignificant on the validation sample but could be characterized more precisely with a larger dataset. While combining these regression, truncation, and discretization uncertainties under an independence assumption would yield a single conservative interval, we do not compute such a joint metric here. This scalar budget addresses only the peak magnitude, meaning spatial localization errors represent a geometric tolerance on the stress zone area that does not affect the scalar comparisons.

Third, our assessment remains restricted to the bounded parameter space. Off-grid evaluation reveals moderate performance degradation, which is more pronounced in the thermal field than the stress field. The physical or numerical mechanism driving this asymmetric sensitivity remains unquantified, requiring a systematic sensitivity analysis of the reduced coordinates across the parameter space in future investigations. Furthermore, the surrogate's accuracy depends on the representativeness of the offline snapshot database, and alternative operating regimes would necessitate constructing an enriched reduced basis.

Finally, the underlying finite-element formulation simplifies the physical brake physics. Neglecting radiative heat transfer, assuming temperature-independent material properties, and prescribing the frictional heat flux rather than solving a fully coupled mechanical contact problem represent major modeling simplifications. The prescribed heat flux assumption likely dominates the overestimation of the peak physical values. This work establishes code verification rather than physical validation against experimental data.

\subsection{Positioning}
\label{subsec:positioning}

The non-intrusive POD-GPR architecture fits within established families of projection-based and data-driven reduced-order models~\cite{benner2015survey}. This framework shares structural features with non-intrusive operator inference and manifold interpolation techniques~\cite{peherstorfer2016operatorinference,guo2019timedependent}, real-time structural surrogate modeling~\cite{mainini2015surrogate}, and non-intrusive reduced-order models incorporating uncertainty quantification~\cite{cicci2023uqromgpr}. The primary contribution of this study lies not in the underlying surrogate architecture, but in the verification methodology that measures emulator errors directly against solver discretization uncertainties.

This methodology distinguishes our framework from adjacent reduced-order modeling paradigms. Certified reduced-basis methods provide mathematically rigorous a posteriori error bounds for field solutions or outputs~\cite{rozza2008rbaposteriori,grepl2005aposteriori}, but they require intrusive access to the solver operators, whereas our non-intrusive budget targets the scalar peak directly against a black-box solver. Bayesian calibration frameworks quantify the discrepancy between numerical predictions and physical experiments~\cite{kennedy2001bayesiancalib}, addressing a surrogate-to-physics gap that remains outside the scope of this work.

Although commercial surrogate modeling toolkits offer comparable spatial field emulation capabilities, our framework explicitly accounts for the approximation error against the solver's own discretization uncertainty. While the computational cost of building the snapshot database scales with the refinement of the high-fidelity reference model, this discretization-anchored verification protocol remains independent of mesh size. The emulator serves as an efficient computational tool, while the verified, floor-anchored error budget provides the mathematical foundation for its engineering deployment.
\section{Conclusion and perspectives}
\label{sec:conclusion}

The primary contribution of this work is the verification methodology itself, rather than the specific emulator architecture. It measures non-intrusive reduced-order emulator accuracy against the numerical discretization uncertainty of the reference solver. By explicitly benchmarking surrogate approximation errors against baseline mesh discretization limits, this approach extends conventional reduced-order model validation procedures. The resulting framework provides a solver-aware fidelity criterion, preventing unnecessary optimization of surrogate models below the numerical resolution already limiting the reference solution. This verification protocol remains applicable to any non-intrusive reduced-order model of parametrized finite-element simulations where a formal discretization error estimate can be established.

We demonstrated this methodology on a carbon-carbon aircraft brake disc, emulating the transient peak temperature and peak von Mises stress fields using proper orthogonal decomposition coupled with Gaussian process regression. Our a posteriori error budget successfully isolates the approximation errors associated with basis truncation from those introduced by coordinate regression. On an independent, blinded test set, the emulator replicates the peak reference solver values within the numerical discretization uncertainty of the mesh for both physical fields, confirming the surrogate-to-code fidelity of the framework. The underlying error budget reveals that linear POD spatial truncation, rather than Gaussian process regression, represents the primary bottleneck to further surrogate refinement.

Two results scope the dimensional reduction for this application. Because the linear POD truncation error falls below the verified discretization floor for the current parameter domain, adopting complex nonlinear reduction schemes, such as deep autoencoders, would provide no observable performance benefits under our verification metric~\cite{greif2019nwidth}. This corollary is specific to the present brake family and parameter domain and does not generalize to all thermomechanical systems. As secondary contributions, the unified non-intrusive workflow predicts both physical fields simultaneously and localizes the peak von Mises stress within a probabilistic spatial zone. Although the exact coordinates of the stress maximum are highly sensitive, the surface area of this high-stress zone varies consistently with the operating parameters, serving as a stable scalar descriptor of stress concentration that the emulator accurately reproduces.

Several numerical and physical boundaries outline the scope of this study. The residual mesh dependence of the peak von Mises stress at the production resolution requires treating the stress discretization floor as a case-dependent range rather than a uniform scalar threshold. Furthermore, the level-one predictive intervals quantify only the epistemic regression uncertainty of the Gaussian processes, omitting spatial truncation and solver discretization errors which our framework handles independently. Consequently, establishing a total predictive variance would require a unified mathematical formulation combining all three error sources. The validated surrogate accuracy remains restricted to the bounded parameter space and exhibits moderate degradation under off-grid extrapolation. Finally, because the high-fidelity finite-element reference assumes temperature-independent material properties and omits radiative heat transfer, this study establishes a code-verification framework rather than a physical validation against experimental measurements.

Future research directions follow directly from these constraints. Resolving the entire offline database on the refined mesh would eliminate residual stress mesh dependencies, converting the case-by-case discretization ranges into a sharp, uniform numerical floor. Integrating the Gaussian process regression covariance, the linear POD truncation bounds, and the discretization error bands into a single, unified uncertainty framework represents an important mathematical challenge that would yield robust confidence intervals for engineering design decisions. Finally, extending this non-intrusive verification methodology from static peak predictions to full spatio-temporal histories requires addressing transient error propagation and time-dependent model reduction, determining whether our surrogate-to-code fidelity criterion remains robust in the presence of temporal error accumulation.

\section*{Data availability}
Upon publication, we will make all supporting data openly available on Zenodo, including the aggregated metrics (in \texttt{numbers.tex}), the five-parameter design of experiments, the exact \nTrain{}/\nVal{}/\nTest{} data splits, and the reproducible software environment. To safeguard proprietary geometric designs, we explicitly exclude the raw four-dimensional thermal and stress tensors, the POD spatial modes, the nodal coordinates, and the underlying computational mesh. This study relies on no proprietary or client-owned datasets. Readers may contact the author to request additional details.

\section*{Declaration of competing interest}
The author works within the aerospace industry. However, this study was pursued strictly as an independent, personal research project. It received no funding, resources, or involvement from the employer, who played no role in designing, executing, analyzing, or reporting this work. The author declares no other competing interests.

\section*{Declaration of generative AI and AI-assisted technologies in the writing process}
The author declares that no scientific content, data analysis, mathematical formulation, or physical interpretation in this manuscript was generated by artificial intelligence. A large language model was used exclusively in a minor capacity as a linguistic aid to polish the phrasing and style of the author's original draft. The author thoroughly reviewed all suggested refinements and assumes sole responsibility for the final publication.

\appendix

\setcounter{figure}{0}
\renewcommand{\thefigure}{A.\arabic{figure}}

\section{Supplementary results}
\label{app:supplementary}

This appendix collects the figures that support the results of the main text without
being essential to the argument. They are organized by theme.

\subsection{Full-field reconstruction}
\label{app:reconstruction}

\begin{figure}[t]
  \centering
  \includegraphics[width=\linewidth]{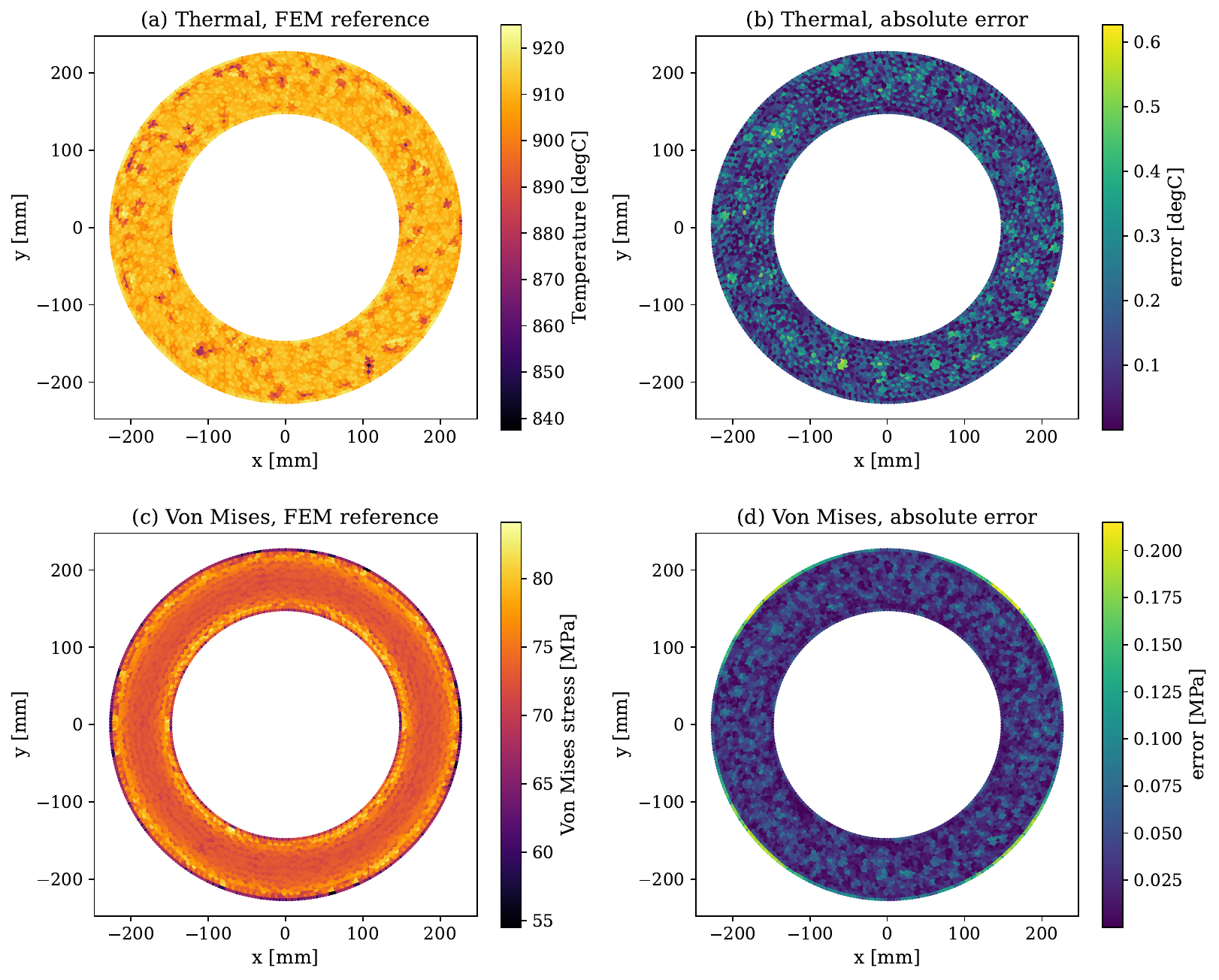}
  \caption{Reference fields and absolute reconstruction errors at the operating rank (thermal field on top, von Mises stress at the bottom). Errors remain concentrated away from the primary friction band. The peak regions are resolved within the limits documented in our main analysis.}
  \label{fig:reconstruction-map}
\end{figure}

\subsection{Predictive-interval calibration}
\label{app:calibration}

\begin{figure}[t]
  \centering
  \includegraphics[width=\linewidth]{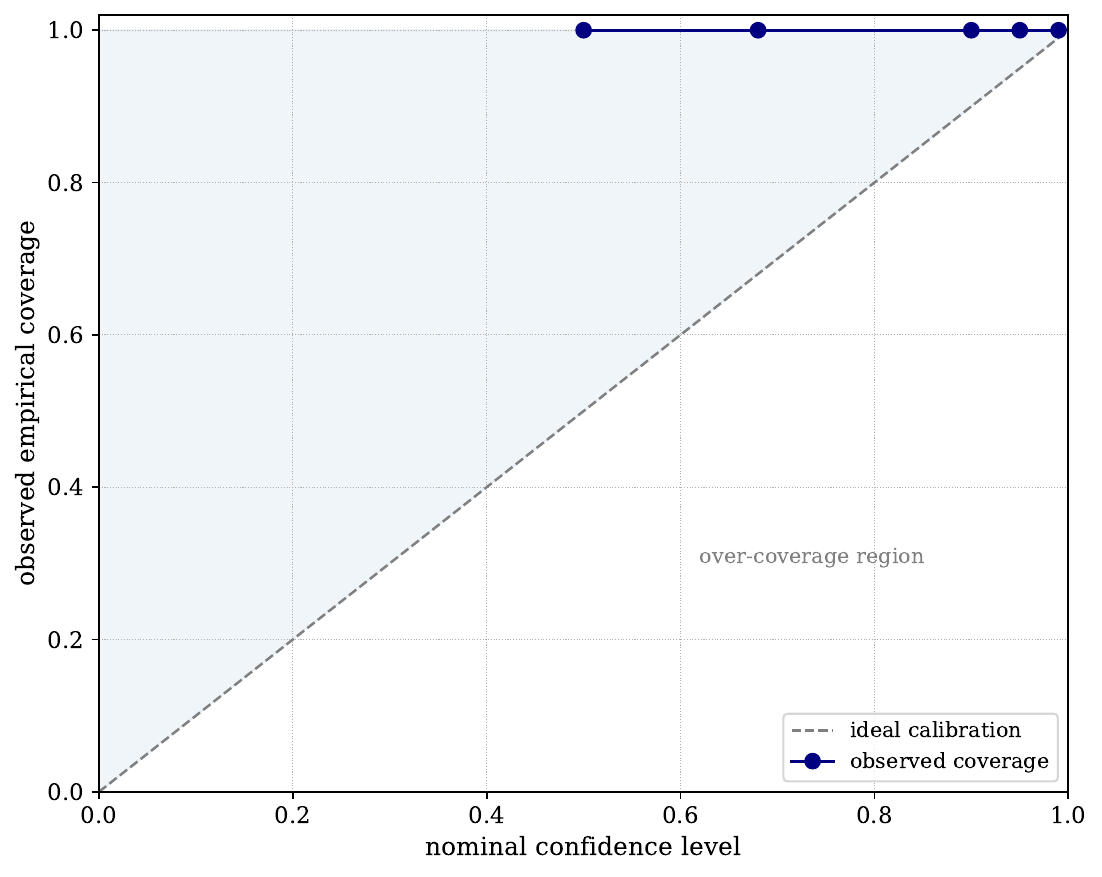}
  \caption{Calibration performance of the level-one predictive interval, comparing actual vs. nominal coverage for both fields using the validation split.}
  \label{fig:interval-calibration}
\end{figure}

\begin{figure}[t]
  \centering
  \includegraphics[width=\linewidth]{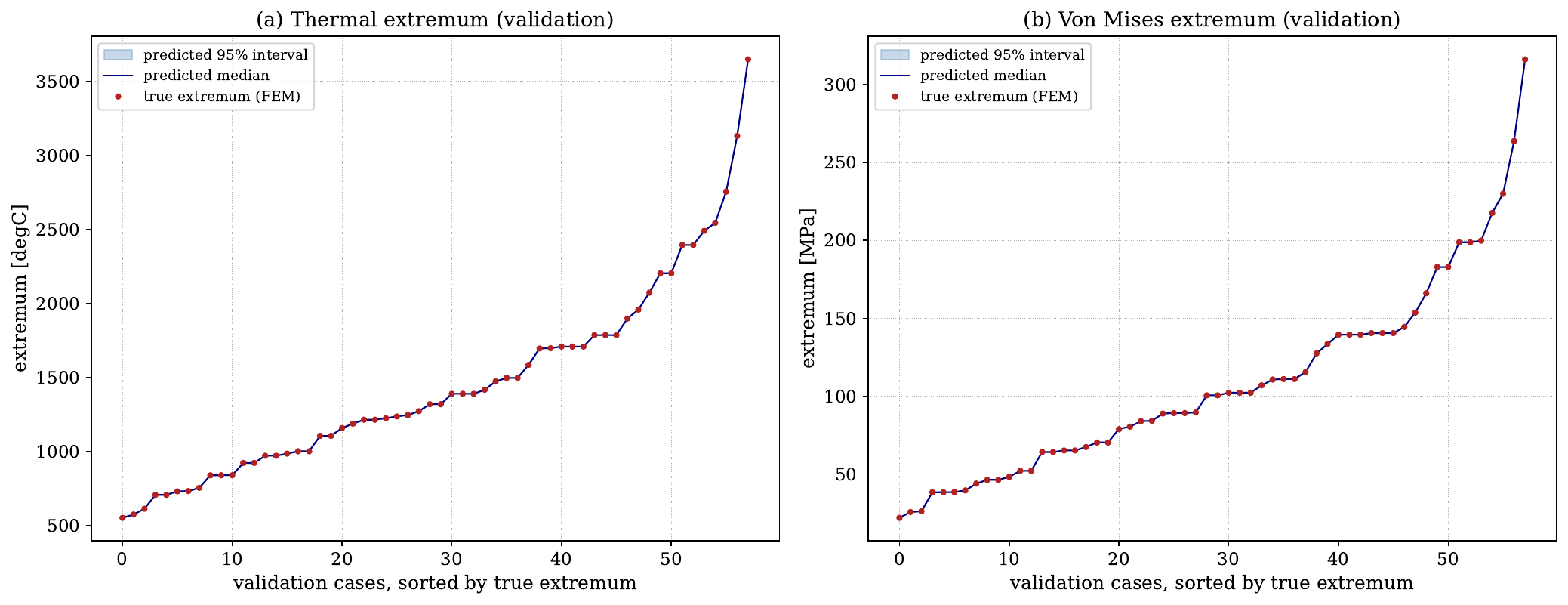}
  \caption{Validation set predictive intervals (level-one) at the extremum. Each case displays the estimated interval alongside the true reference peak.}
  \label{fig:intervals-validation}
\end{figure}

\begin{figure}[t]
  \centering
  \includegraphics[width=\linewidth]{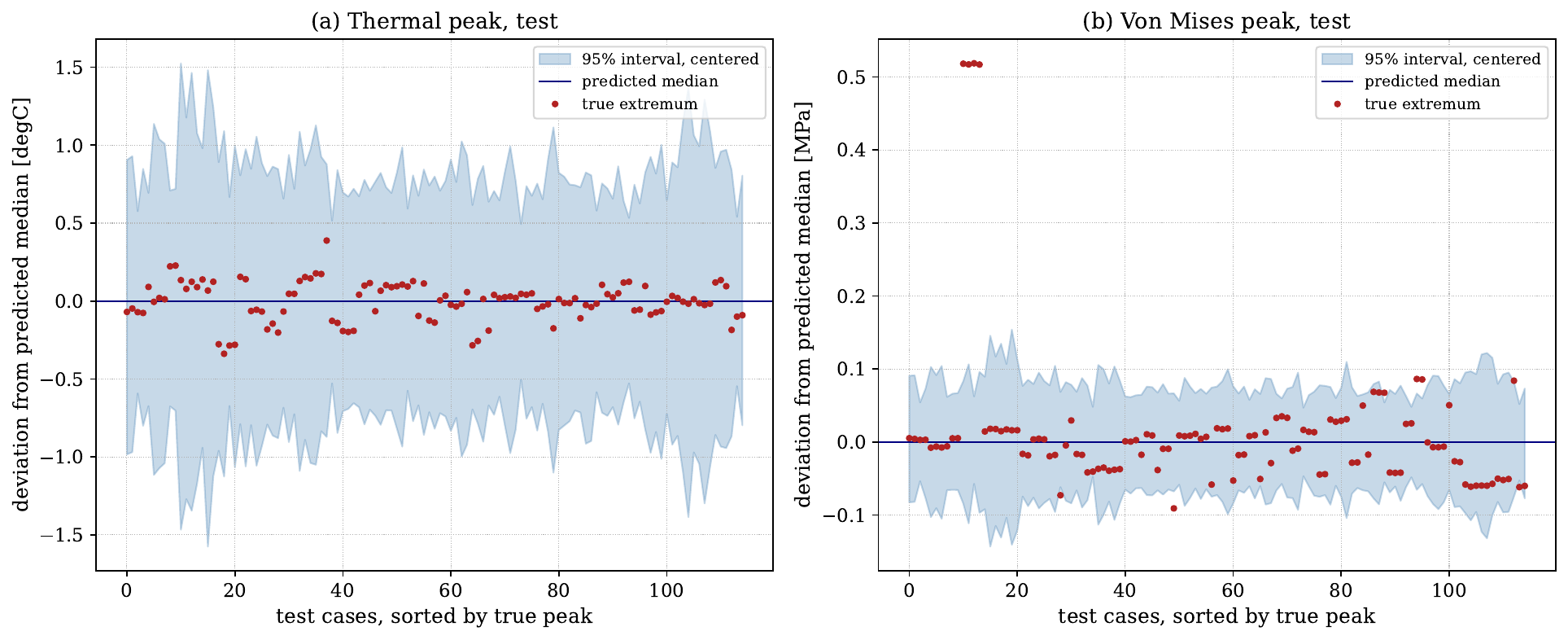}
  \caption{Predictive intervals (level-one) at the extremum on the independent test set. True peaks are plotted against their corresponding predicted intervals.}
  \label{fig:intervals-test}
\end{figure}

\subsection{Extremum error}
\label{app:extremum-error}

\begin{figure}[t]
  \centering
  \includegraphics[width=\linewidth]{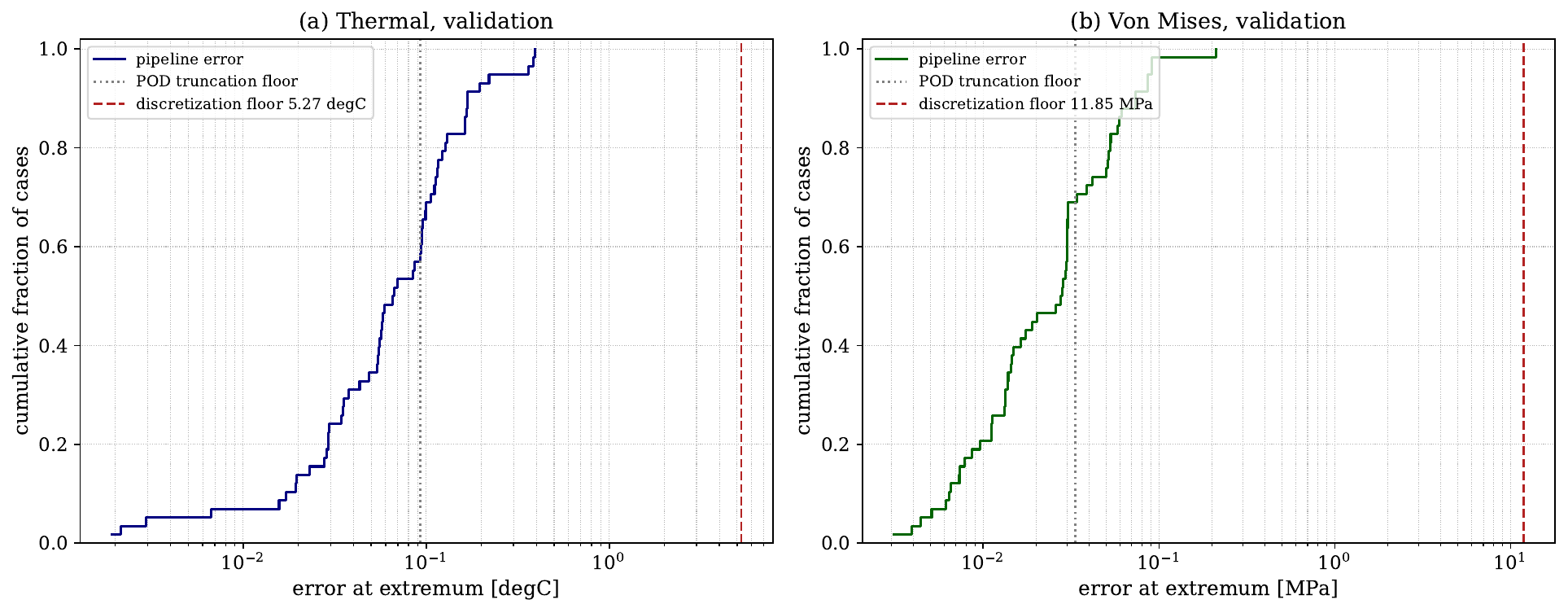}
  \caption{Relative extremum errors across the validation set, sorted by case. The horizontal line indicates our pre-registered conformity threshold.}
  \label{fig:extremum-error-validation}
\end{figure}

\subsection{Computational speedup}
\label{app:speedup}

\begin{figure}[t]
  \centering
  \includegraphics[width=\linewidth]{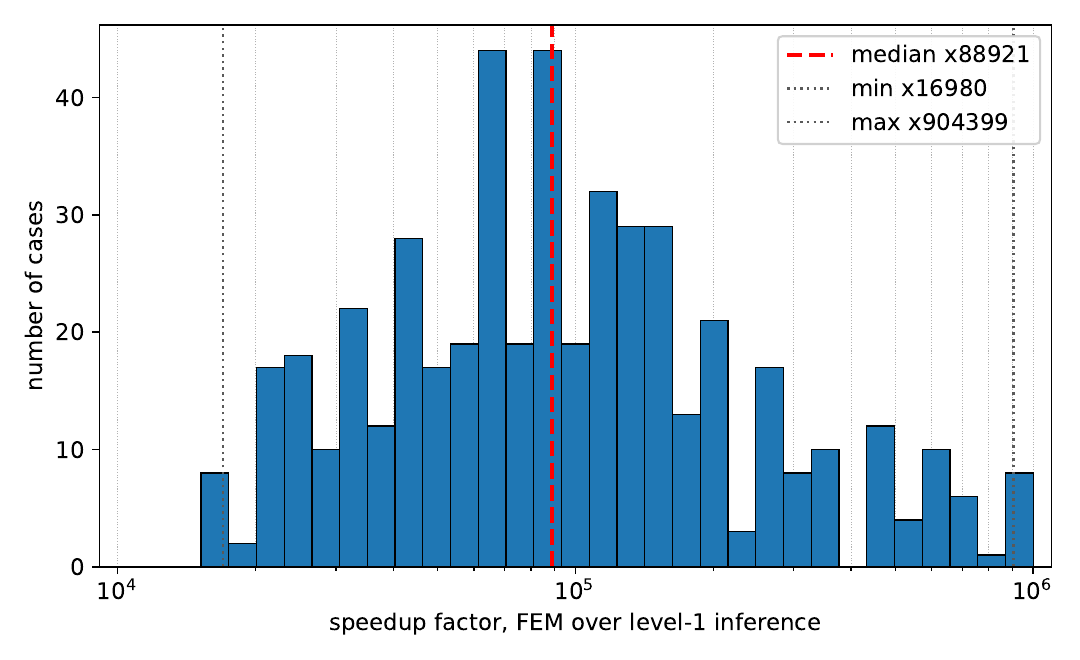}
  \caption{Speedup factor distribution of our emulator relative to the finite-element solver (measured on the reference hardware, excluding disk I/O latency).}
  \label{fig:speedup-distribution}
\end{figure}

\begin{figure}[t]
  \centering
  \includegraphics[width=\linewidth]{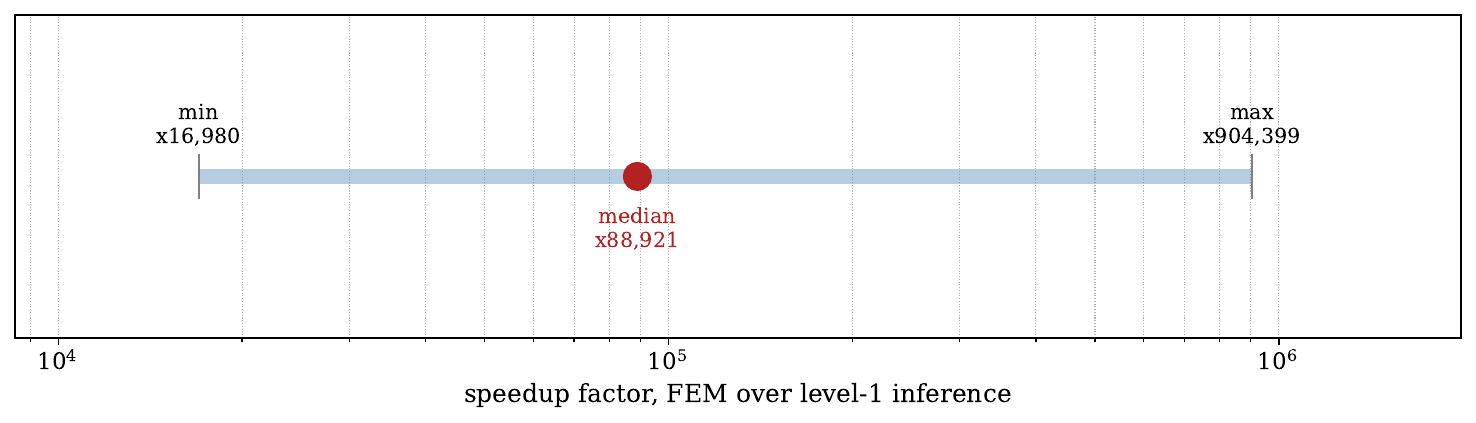}
  \caption{Emulator speedup compared to solver computational cost on the reference hardware. Rather than a flat metric, speedup varies dynamically because physical braking durations scale the baseline solver cost.}
  \label{fig:speedup-range}
\end{figure}

\subsection{Residual correlation and the independence assumption}
\label{app:residual}

\begin{figure}[t]
  \centering
  \includegraphics[width=\linewidth]{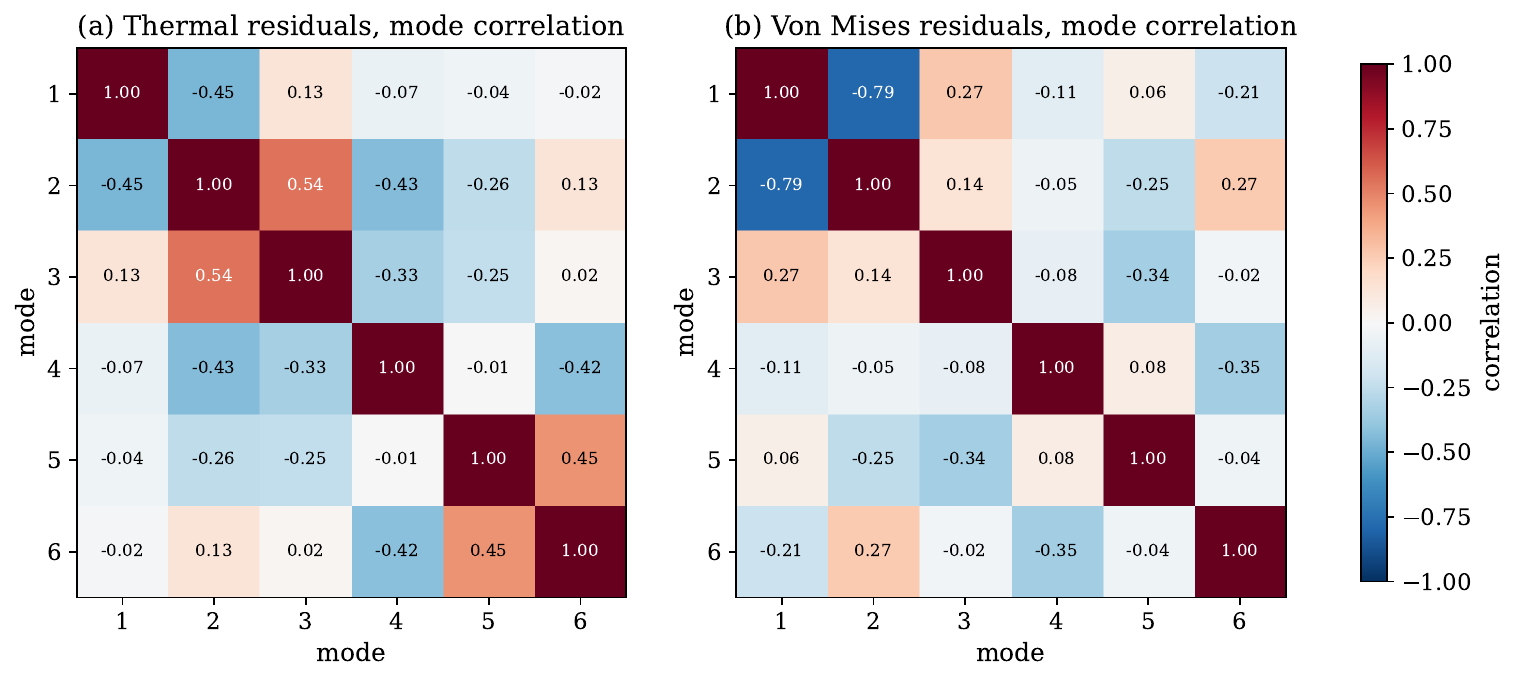}
  \caption{Cross-correlation of the validation residuals across off-diagonal reduced dimensions for the von Mises field.}
  \label{fig:residual-correlation}
\end{figure}

\begin{figure}[t]
  \centering
  \includegraphics[width=\linewidth]{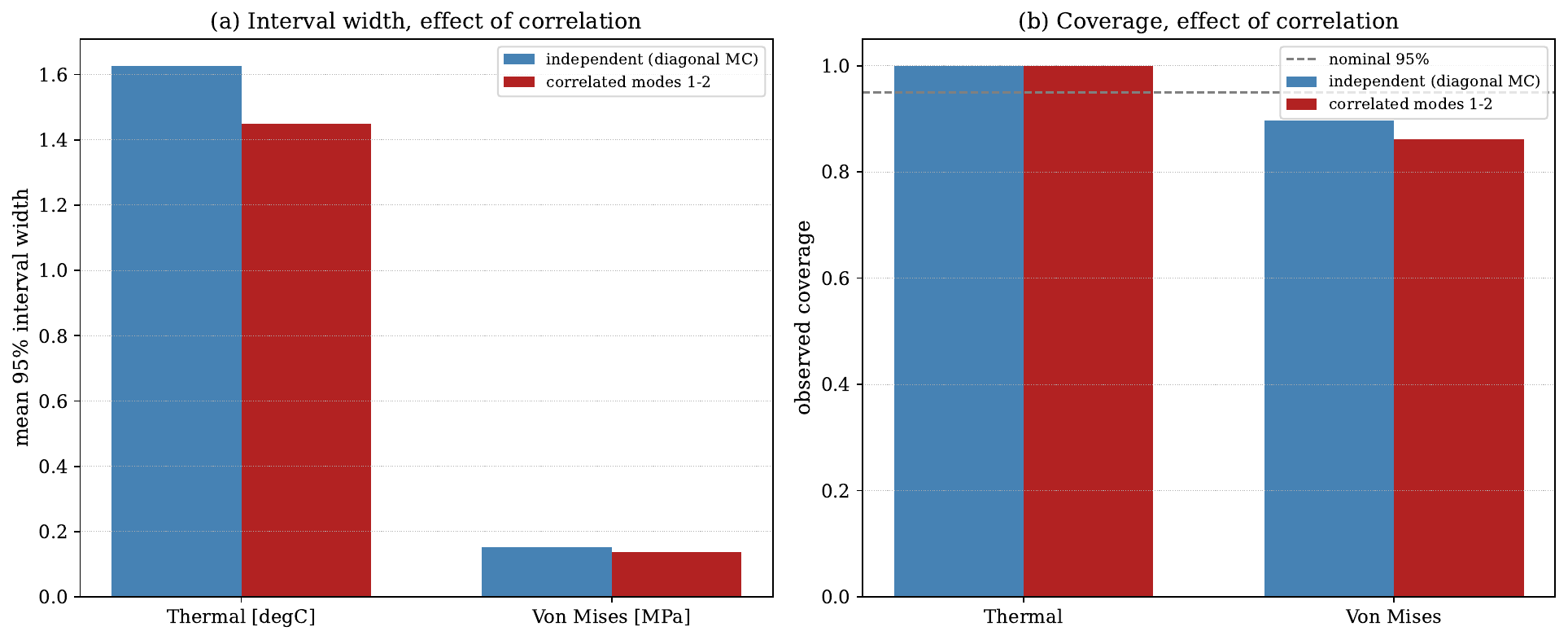}
  \caption{Influence of the coordinate independence assumption on extremum coverage. We compare a diagonal Monte Carlo approach with full-covariance joint sampling.}
  \label{fig:correlation-impact}
\end{figure}

\begin{figure}[t]
  \centering
  \includegraphics[width=\linewidth]{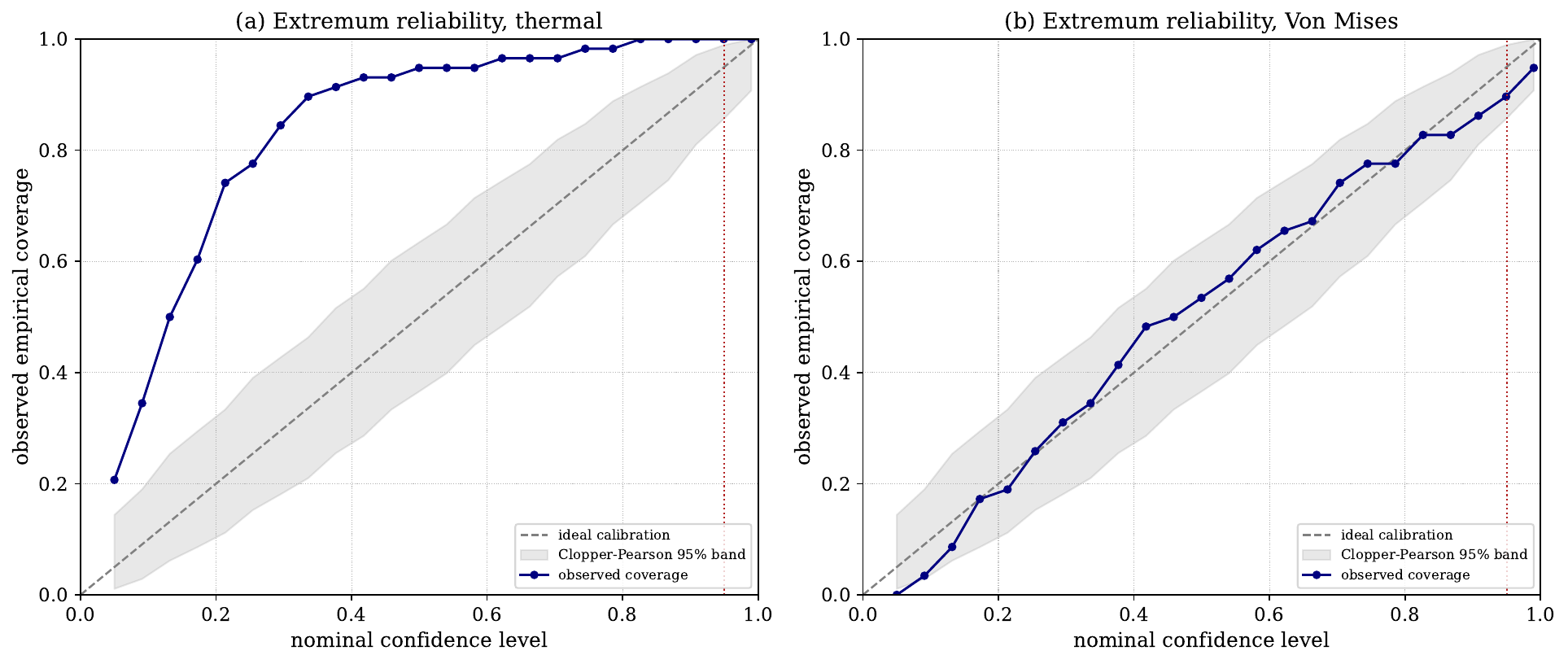}
  \caption{P-P plot showing the reliability of the level-one peak interval, matching observed coverage directly to nominal expectations.}
  \label{fig:extremum-reliability}
\end{figure}

\subsection{Partition integrity}
\label{app:partition}

\begin{figure}[t]
  \centering
  \includegraphics[width=\linewidth]{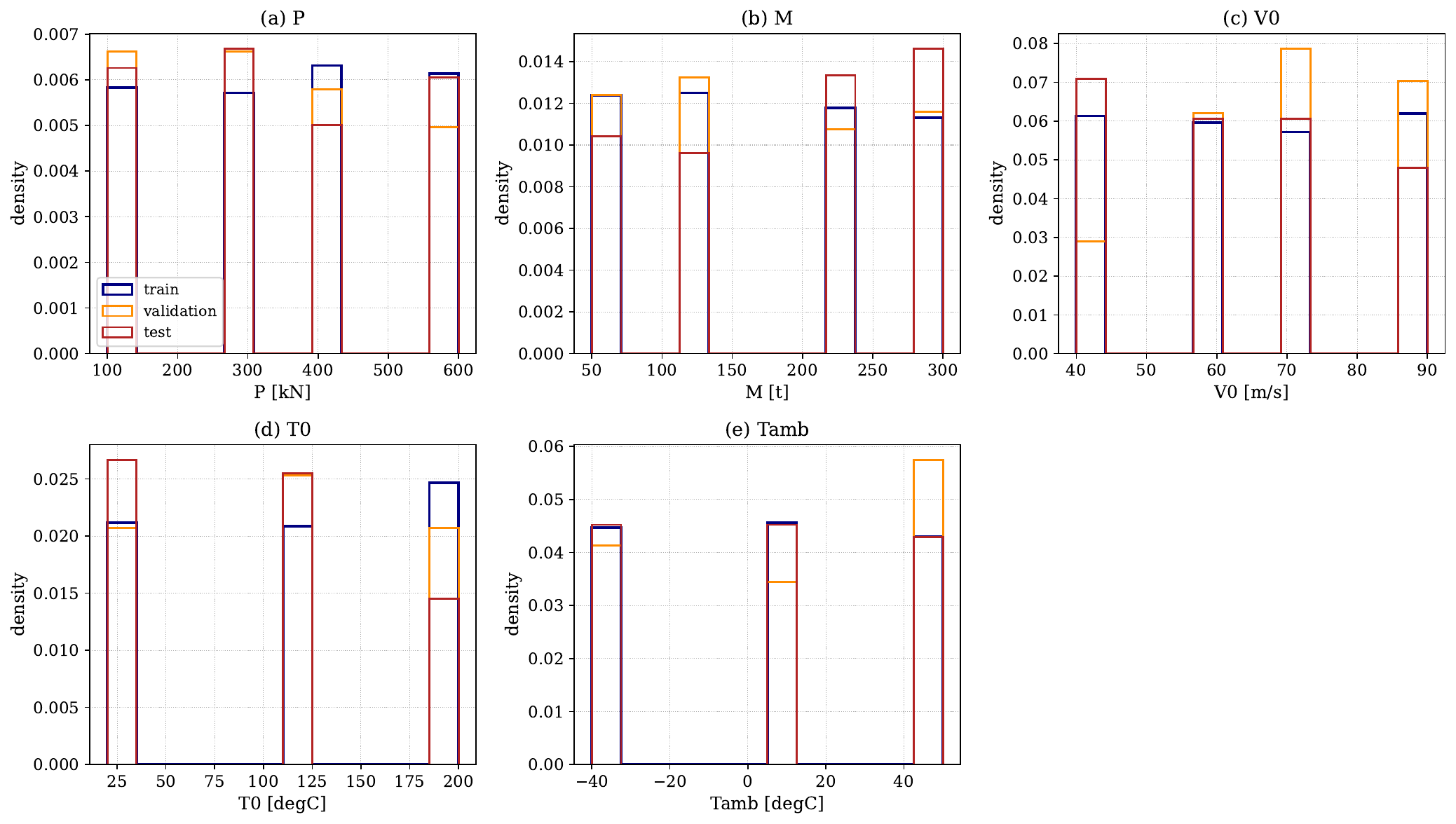}
  \caption{Parameter level frequencies across our three data splits, verifying that the partition maintains design-of-experiment balance without any data leakage.}
  \label{fig:partition-no-leakage}
\end{figure}

\subsection{Zonal localization}
\label{app:zone}

\begin{figure}[t]
  \centering
  \includegraphics[width=\linewidth]{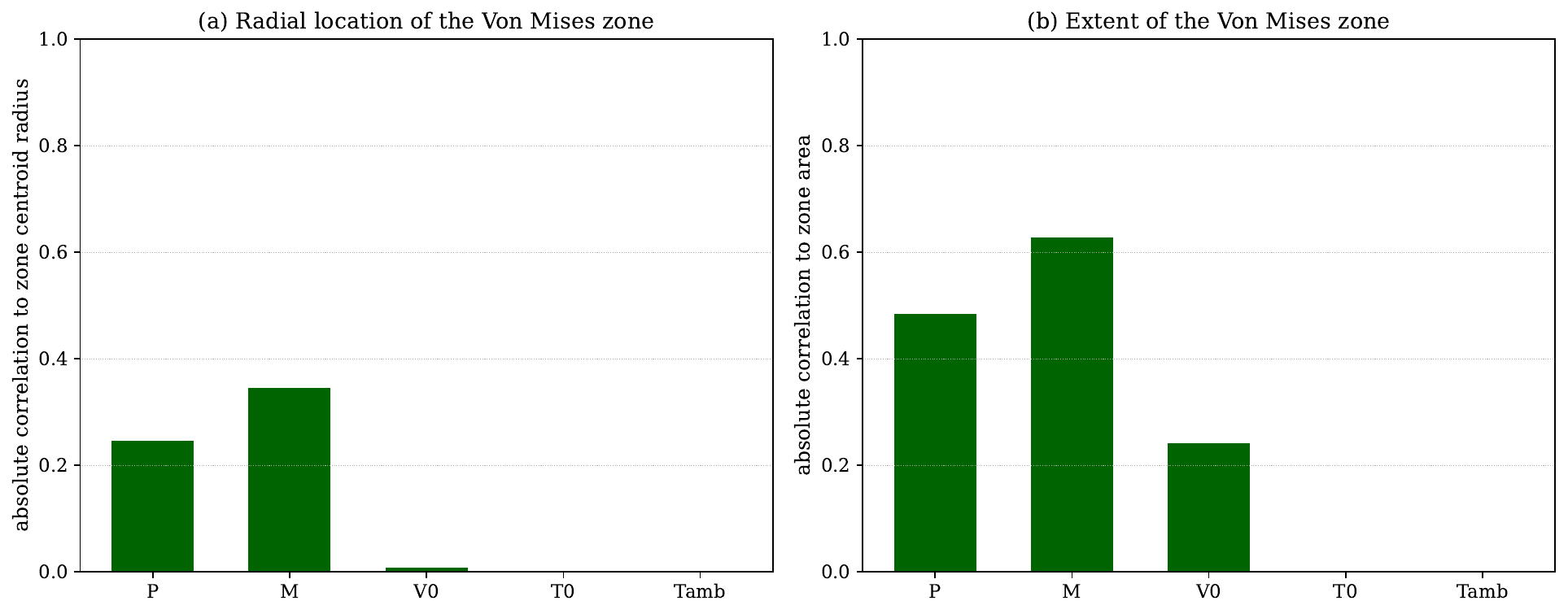}
  \caption{How operating conditions scale the predicted von Mises zone area. The area changes alongside mass, braking load, and initial speed, remaining insensitive to ambient and starting temperatures.}
  \label{fig:zone-variation}
\end{figure}

\begin{figure}[t]
  \centering
  \includegraphics[width=\linewidth]{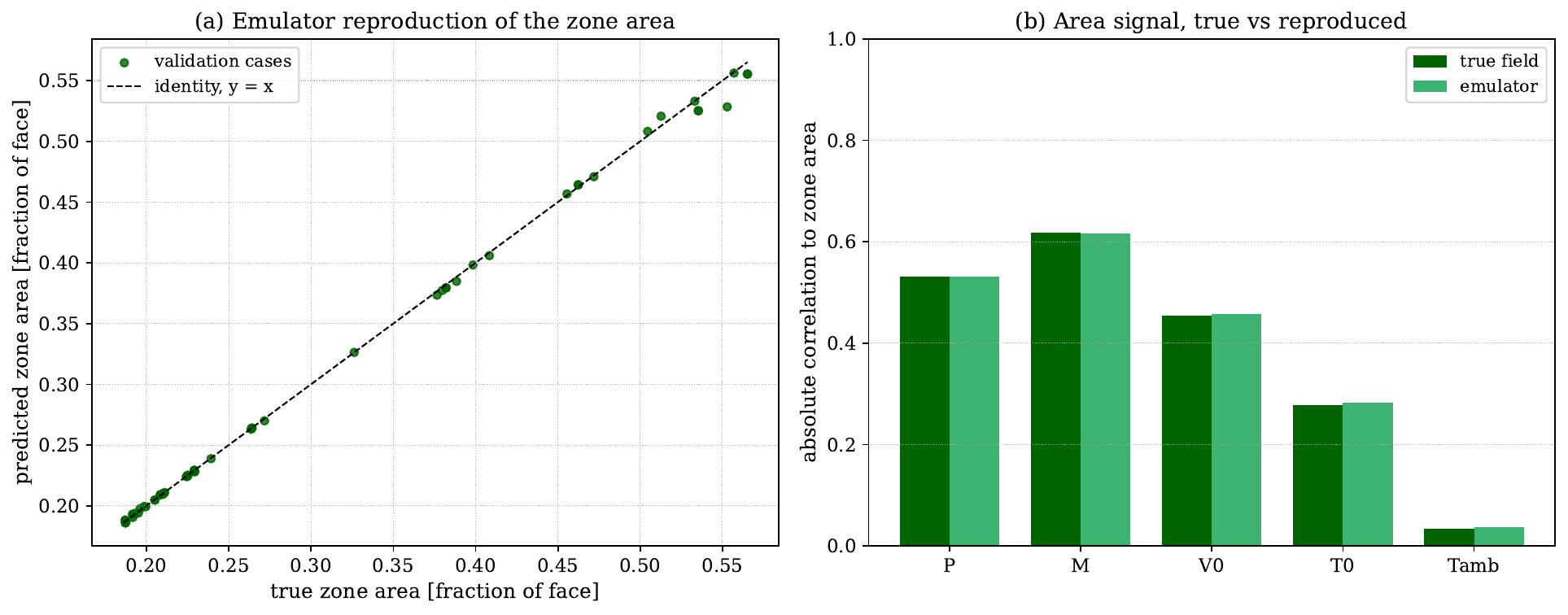}
  \caption{Validation set comparison of predicted vs. reference von Mises zone areas, demonstrating reconstruction within our defined target tolerance.}
  \label{fig:area-reproduction}
\end{figure}

\begin{figure}[t]
  \centering
  \includegraphics[width=\linewidth]{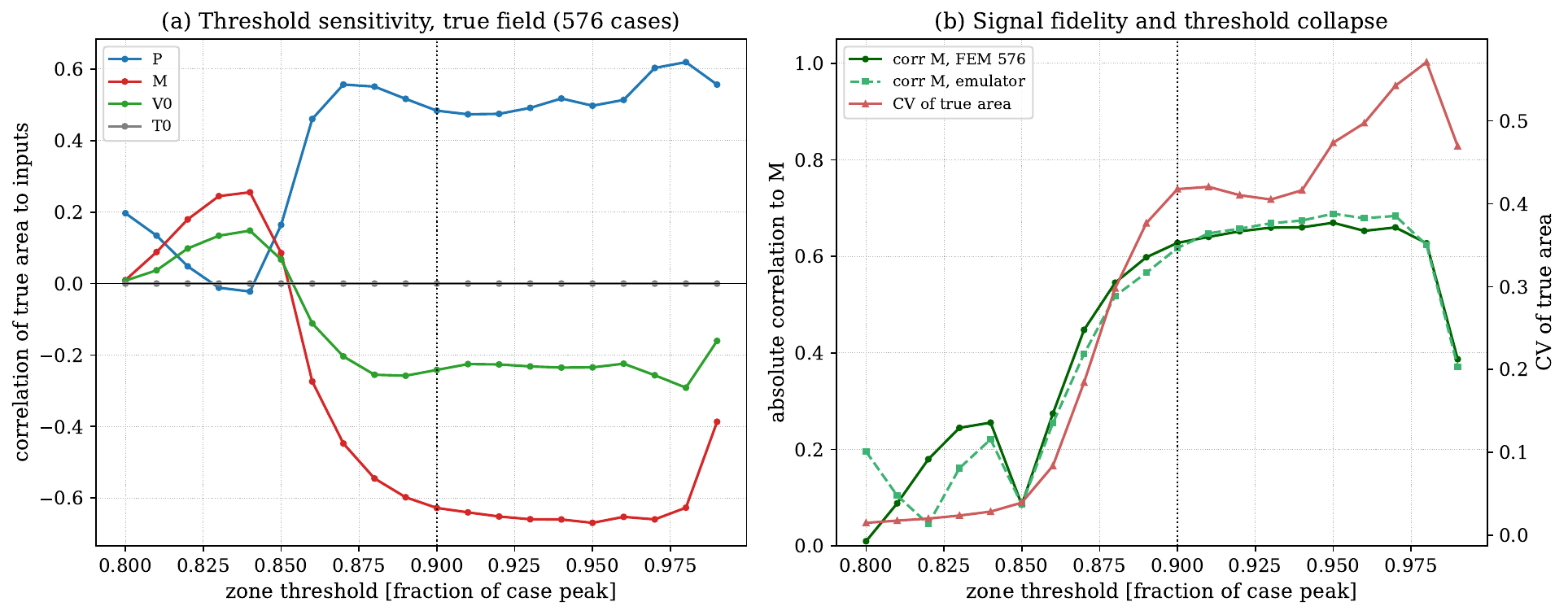}
  \caption{Sensitivity analysis of zone shape metrics across varying super-level-set thresholds, showing stable area metrics and parameter trends.}
  \label{fig:threshold-robustness}
\end{figure}

\clearpage

\bibliographystyle{elsarticle-num}
\bibliography{references}

\end{document}